\documentclass[12pt]{amsart}

\newtheorem{theorem}{Theorem}[section]
\newtheorem{lemma}[theorem]{Lemma}
\newtheorem{corollary}[theorem]{Corollary}
\newtheorem{proposition}[theorem]{Proposition}
\newtheorem{remark}[theorem]{Remark}
\newtheorem{definition}[theorem]{Definition}

\newcommand{\ncom}{\newcommand}
\ncom{\lrar}{\longrightarrow}
\ncom{\ov}{\overline}
\ncom{\m}{\mbox}
\ncom{\sta}{\stackrel}
\ncom{\comx}{{\mathbb C}}
\ncom{\Z}{{\mathbb Z}}
\ncom{\Q}{{\mathbb Q}}
\ncom{\R}{{\mathbb R}}
\ncom{\G}{{\mathbb G}}
\ncom{\al}{\alpha}
\ncom{\p}{{\mathbb P}}
\ncom{\E}{{\mathbb E}}
\ncom{\N}{{\mathbb N}}
\ncom{\K}{{\mathbb K}}
\ncom{\bbH}{{\mathbb H}}
\ncom{\bbL}{{\mathbb L}}
\ncom{\f}{\frac}
\ncom{\cA}{{\mathcal A}}
\ncom{\cB}{{\mathcal B}}
\ncom{\cD}{{\mathcal D}}
\ncom{\cX}{{\mathcal X}}
\ncom{\cO}{{\mathcal O}}
\ncom{\cW}{{\mathcal W}}
\ncom{\cL}{{\mathcal L}}
\ncom{\cP}{{\mathcal P}}
\ncom{\cH}{{\mathcal H}}
\ncom{\cS}{{\mathcal S}}
\ncom{\cM}{{\mathcal M}}
\ncom{\cC}{{\mathcal C}}
\ncom{\cT}{{\mathcal T}}
\ncom{\cF}{{\mathcal F}}
\ncom{\cN}{{\mathcal N}}
\ncom{\cJ}{{\mathcal J}}
\ncom{\cV}{{\mathcal V}}
\ncom{\cZ}{{\mathcal Z}}
\ncom{\cU}{{\mathcal U}}
\ncom{\cSU}{{\mathcal S \mathcal U}}
\ncom{\cG}{{\mathcal G}}
\ncom{\cQ}{{\mathcal Q}}
\ncom{\cR}{{\mathcal R}}
\ncom{\cY}{{\mathcal Y}}

\ncom{\eop}{{\hfill $\Box$}}

\begin{document}
\baselineskip=16pt

\title[Parabolic Chern character of the de Rham bundles]{A relation between the parabolic Chern characters of the de Rham bundles}
\author[J. N. Iyer]{Jaya NN. Iyer}
\address{The Institute of Mathematical Sciences, CIT
Campus, Taramani, Chennai 600113, India}
\email{jniyer@imsc.res.in}

\author[C. T.  Simpson]{Carlos T.  Simpson}
\address{CNRS, Laboratoire J.-A.Dieudonn\'e, Universit\'e de Nice--Sophia Antipolis,
Parc Valrose, 06108 Nice Cedex 02, France}
\email{carlos@math.unice.fr}

\footnotetext{Mathematics Classification Number: 14C25, 14D05, 14D20, 14D21 }
\footnotetext{Keywords: Connections, Chow groups, de Rham cohomology, parabolic bundles, Lefschetz pencils.}
\begin{abstract}
In this paper, we consider the weight $i$ de Rham--Gauss--Manin bundles 
on a smooth variety arising from a smooth 
projective morphism $f:X_U\lrar U$ for $i\geq 0$. We associate to each weight $i$ 
de Rham bundle, a certain parabolic bundle on $S$
and consider their parabolic Chern characters in the rational Chow groups, for a 
good compactification $S$ of $U$. We
show the triviality of the alternating sum of these parabolic bundles in the 
(positive degree) rational Chow groups. This removes the hypothesis of
semistable reduction in the original result of this kind due to Esnault and Viehweg. 
\end{abstract}

\maketitle

$$\textsc{Contents}$$

1. Introduction

2. Parabolic bundles

3. The parabolic bundle associated to a logarithmic connection

4. Lefschetz fibrations

5. Main Theorem

6. Appendix: an analogue of Steenbrink's theorem

7. References

\section{Introduction}
Suppose $X$ and $S$ are  irreducible projective varieties defined over the complex numbers and
$\pi:X\lrar S$ is a morphism  such that the restriction
$X_U\rightarrow U$ over a nonsingular dense open set is smooth of relative dimension $n$. 
The following bundle on $U$, for $i\geq 0$,
$$
\cH^i:=R^i\pi_*\Omega^\bullet_{X_U/U}
$$
is equipped with a flat connection $\nabla$, called as the Gauss-Manin connection.
We call the pair $(\cH^i,\nabla)$ the de Rham bundle or the Gauss-Manin bundle of weight $i$.

Suppose $S$ is a nonsingular compactification of $U$ such that $D:=S -U$ is a normal crossing 
divisor and the associated local system $\m{ker}\nabla$ has unipotent monodromies along the components of $D$.
The bundle $\ov\cH^i$
is the \textit{canonical extension} of $\cH^i$ (\cite{De}) and is equipped with a logarithmic flat 
connection $\ov\nabla$. It is characterised by the property that it has nilpotent residues.

By the Chern--Weil theory, the de Rham Chern classes
$$
c^{dR}_i(\cH^k)\,\in\, H^{2i}_{dR}( U)
$$
vanish, and by a computation shown in
\cite[Appendix B]{Es-Vi1}, the de Rham classes
$$
c^{dR}_i(\overline{\cH}^k) \,\in\, H^{2i}_{dR}(S)
$$
vanish too. The essential fact used is that the residues of
$\overline{\nabla}$ are nilpotent.

The algebraic Chern--Simons theory initiated by S. Bloch and H. Esnault (\cite{BE})
studies the Chern classes (denoted by $c^{Ch}_i$) of flat bundles in
the rational
Chow groups of $U$ and $S$. It is conjectured by H. Esnault
that the classes $c^{Ch}_i(\cH^{k})$ and
$c^{Ch}_i(\overline{\cH}^{k})$ are trivial
for all $i>0$ and $k$ (\cite[p. 187--188]{Es1}, \cite{Es2}).

The cases where it is known to be true are as follows.
In \cite{Mu}, Mumford proved this for any family of stable curves.
In \cite{vdG}, van der Geer proved that
$c^{Ch}_i(\cH^{1})$ is trivial
when $X\, \longrightarrow\, S$ is a family of abelian
varieties. Further, for any family of abelian varieties
of dimension $g$, the
rational Chow group elements $c^{Ch}_i(\overline{\cH}^{1})$,
$i\geq 1$, were proved to be trivial
by Iyer under the assumption that $g\, \leq \, 5$ (\cite{Iy})
and by Esnault and Viehweg and for all $g\,>\,0$ (\cite{Es-Vi2}).
Further, for some families of moduli spaces, Biswas and Iyer (\cite{Bi-Iy}) have checked the triviality
of the classes in the rational Chow groups.

In this paper we consider parabolic bundles associated to logarithmic connections (section \ref{Parlog}) 
instead of canonical extensions.
By Steenbrink's theorem \cite[Proposition 2.20]{Steen} (see also \cite{KatzICM}), the monodromy of the VHS associated to families 
is quasi-unipotent and the residues have rational eigenvalues. It is natural to consider the parabolic 
bundles associated to such local systems. Further, these parabolic bundles are compatible with 
pullback morphisms (Lemma \ref{pullback}), unlike canonical extensions.
One can also define the Chern character of parabolic bundles in the rational Chow groups (section \ref{ChernPar}). 
All this is possible
by using a correspondence of these special parabolic bundles, termed as locally abelian parabolic bundles, 
with  vector bundles on a particular DM-stack (Lemma \ref{equivParaDM}).

In this framework, we show

\begin{theorem}\label{th.-1}
Suppose $\pi:X_U\lrar  U$ is a smooth
projective morphism of relative dimension $n$ between
nonsingular varieties. Consider a nonsingular compactification
$U\subset S$ such that $S -U$ is a normal crossing divisor.
Then the Chern character of the alternating sum of the parabolic bundles $\ov\cH^i(X_U/U)$ in each degree,
$$
\sum _{i=0}^{2n} (-1)^i{\rm ch}(\overline{\cH}^i(X_U/U))
$$
lies in $CH^0(S)_\Q$ or equivalently the pieces in all of the positive-codimension Chow groups with rational coefficients
vanish.

Here $\ov\cH^i(X_U/U)$ denotes the parabolic bundle associated to the weight $i$ de Rham bundle on $U$.
\end{theorem}

In fact we will prove the same thing when the morphism is not generically smooth, where $U$ is the open set
over which the map is topologically a fibration, see Theorem \ref{th.-deRham}. 

If $X_U\rightarrow U$ has a semi-stable extension $X\rightarrow S$ (or a compactified family satisfying certain conditions) 
then the triviality of the Chern character of the alternating sum of de Rham bundles
in the (positive degree) rational Chow groups is proved by Esnault and Viehweg in \cite[Theorem 4.1]{Es-Vi2}. 
This is termed as a logarithmic Grothendieck
Riemann-Roch theorem (GRR) since GRR is applied to the logarithmic relative de Rham sheaves to obtain the relations.
It might be possible to generalize the calculation of \cite{Es-Vi2} to the case of a weak toroidal semistable
reduction which always exists by \cite{AbramovichKaru}. However, this seems like it would be difficult to set up. 

We generalise the Esnault-Viehweg result and a good compactified family is not required over $S$. In particular, 
we do not use calculations with the Grothendieck Riemann-Roch formula, although we do use the
general existence of such a formula. Further, Theorem \ref{th.-1} shows that
the singularities in the fibres of extended families do not play any role, in higher dimensions. Instead, our inductive
argument calls upon de Jong's semistable reduction for curves \cite{deJong} at each inductive step. 

As an application, we show 
\begin{corollary}\label{th.-2}
Suppose $X_U\lrar U$ is any 
family of projective surfaces and $S$ is a good compactification of $U$. 
Then the parabolic Chern character satisfies
$$ch(\ov\cH^i(X_U/U)) \in CH^{0}(S)_\Q
$$
for each $i\geq 0$.
\end{corollary}

This is proved in \S{5.5} Proposition \ref{pr.-surfaces}. On a non-compact base $S$ supporting a smooth family of surfaces, 
this was observed in \cite[Example 7.3]{BE}.
We use the weight filtration on the cohomology of the singular surface, a 
resolution of singularities, the triviality of the classes of $\ov\cH^1$ (\cite{Es-Vi2}),
and using Theorem \ref{th.-1} we deduce the proof.

The proof of Theorem \ref{th.-1} is by induction principle and using the Lefschetz theory (such an approach was used
earlier in \cite{BE2} for a similar question).
We induct on the relative dimension $n$. By the Lefschetz theory, the cohomology of an $n$-dimensional nonsingular
projective variety is expressed in terms of the cohomology of a nonsingular hyperplane section and the cohomology
of the (extension of) \textit{variable local system} on $\p^1$.  This  helps us to apply induction and conclude 
the relations between the Chern
classes of the de Rham bundles in the rational Chow groups (Theorem \ref{th.-1}).
Our proof requires a certain amount of machinery such as the notion of parabolic bundle. The reason for this is
that the local monodromy transformations of a Lefschetz pencil whose fiber dimension is even, are reflections of
order two rather than unipotent transvections. In spite of this machinery we feel that the proof is basically
pretty elementary, and in particular it doesn't require us to follow any complicated calculations with GRR.

A result of independent interest is Theorem \ref{th.-steen} in
\S 6, Appendix, which was written by the second author but was 
occasioned by the talk in Nice given by the first author. This is an
analogue of Steenbrink's theorem \cite[Theorem 2.18]{Steen}. In
this case, the relative dimension is one and the relative
(logarithmic) de Rham complex has coefficients in an extension of
a unipotent local system. It is proved that the associated
cohomology sheaves are locally free and having a Gauss--Manin
connection. Furthermore, the logarithmic connection has nilpotent
residues along the divisor components (where it has poles).

{\small \textit{Acknowledgements}
The first named author would like to thank H. Esnault for introducing the questions on de Rham bundles to her and 
for having useful
conversations at different periods of time. She also thanks A. Hirschowitz for the invitation to visit Nice during 
Dec.2004, when this work and collaboration was begun. The visit was funded by NBHM and CNRS. Both authors again thank H. Esnault for numerous helpful comments correcting errors in the
first version, and pointing out the
reference to \cite{KatzICM} which clarifies the discussion in \S 6. }

\section{Parabolic bundles}

We treat some preliminaries on the notion of parabolic bundle \cite{Seshadri}. 
This takes a certain amount of space, and we are
leaving without proof many details of the argument. The purpose of this discussion in our proof of the main theorem
is to be able to treat the case of Lefschetz pencils of
even fiber dimension, in which case the monodromy transformations are reflections of order two (rather than
the more classical unipotent transformations in the case of odd fiber dimension). Thus we will at the end 
be considering parabolic structures with weights $0, \frac{1}{2}$ and the piece of weight $\frac{1}{2}$ will
have rank one. Furthermore we will assume by semistable reduction that the components of the divisor
of singularities don't touch each other. Nonetheless, it seems better to give a sketch of the general theory
so that the argument can be fit into a proper context.

Suppose $X$ is a smooth variety and
$D$ is a normal crossings divisor. Write $D=\bigcup _{i=1}^k D_i$ as a union of irreducible components,
and we assume that the $D_i$ are themselves smooth meeting transversally.

We will define the notion of {\em locally abelian parabolic bundle on
$(X,D)$}. We claim that this is the right definition of this notion; however intermediate definitions may or may
not be useful (e.g. our notion of ``parabolic sheaf'' might not be the right one).
The notion which we use here appeared for example in Mochizuki \cite{Mochizuki}, and is slightly different from the
one used by Maruyama and Yokogawa \cite{MaruyamaYokogawa} in that we consider different filtrations for all the different
components of the divisor. 

Also we shall only consider parabolic structures with rational weights.

A {\em parabolic sheaf} on $(X,D)$ is a collection of torsion-free sheaves $F_{\alpha}$  indexed by
multi-indices $\alpha =(\alpha _1,\ldots , \alpha _k)$ with $\alpha _i\in \Q$, together with
inclusions of sheaves of $\cO _X$-modules
$$
F_{\alpha} \hookrightarrow F_{\beta}
$$
whenever $\alpha _i\leq \beta _i$ (a condition which we write as $\alpha \leq \beta$ in what follows),
subject to the following hypotheses:
\newline
---(normalization/support)\, let $\delta ^i$ denote the multiindex
$\delta ^i_i=1,\;\; \delta ^i_j=0, \; i\neq j$,
then $F_{\alpha + \delta ^i} = F_{\alpha } (D_i)$ (compatibly with the inclusion); and
\newline
---(semicontinuity)\, for any given $\alpha$ there exists $c>0$ such that for any multiindex $\varepsilon$
with $0\leq \varepsilon _i < c$ we have $F_{\alpha + \varepsilon} = F_{\alpha}$.

It follows from the normalization/support condition that the quotient sheaves $F_{\alpha} / F_{\beta} $
for $\beta \leq \alpha$  are supported in a schematic neighborhood of the divisor $D$, and indeed if
$\beta \leq \alpha \leq \beta + \sum n_i\delta ^i$ then $F_{\alpha} / F_{\beta}$ is supported over
the scheme $\sum _{i=1}^k n_iD_i$. Let $\delta  := \sum _{i=1}^k\delta ^i$. Then
$$
F_{\alpha - \delta} = F_{\alpha} (-D)
$$
and $F_{\alpha} / F_{\alpha -\delta} = F_{\alpha} |_D$.

The semicontinuity condition means that the structure is determined by the sheaves $F_{\alpha}$ for
a finite collection of indices $\alpha$ with $0\leq \alpha _i < 1$, the {\em weights}.

For each component $D_i$ of the divisor $D$, we have
$$
F_{\alpha} |_{D_i} = F_{\alpha} / F_{\alpha -\delta ^i}.
$$
Thus for $\alpha = 0$ we have
$$
F_0|_{D_i} = F_0 / F_{-\delta ^i}.
$$
This sheaf over $D_i$ has a filtration by subsheaves which are the images of the $F_{\beta _i\delta ^i}$ for
$-1 < \beta _i \leq 0$. The filtration stops at $\beta _i=-1$ where by definition the subsheaf is zero.
Call the image subsheaf $F_{D_i, \beta _i}$ and for $\beta _i > -1$ put
$$
Gr _{D_i,\beta _i}:= F_{D_i, \beta _i}/ F_{D_i, \beta _i-\varepsilon }
$$
with $\varepsilon$ small. There are finitely many values of $\beta _i$ such that the $Gr$ is nonzero. These
values are the {\em weights} of $F$ along the component $D_i$. Any $\Z$-translate of one of these weights
will also be called a weight. The global parabolic structure is determined by the sheaves $F_{\alpha}$ for
multiindices $\alpha$ such that each $\alpha_i$ is a weight along the corresponding component $D_i$.

\subsection{The locally abelian condition}\label{lap}

A {\em parabolic line bundle} is a parabolic sheaf $F$ such that all the $F_{\alpha}$ are line bundles.
An important class of examples is obtained as follows: if $\alpha$ is a multiindex then we can define
a parabolic line bundle denoted
$$
F:= \cO _X(\sum _{i=1}^k \alpha _i D_i)
$$
by setting
$$
F_{\beta } := \cO _X (\sum _{i=1}^k a_i D_i)
$$
where each $a_i$ is the largest integer such that $a_i\leq \alpha _i + \beta _i$.

If $F$ is a parabolic sheaf, set $F_{\infty}$ equal to the extension $j_{\ast}(j^{\ast} F_{\alpha})$ for
any $\alpha$, where $j:X-D\hookrightarrow X$ is the inclusion. It is the sheaf of sections of $F$ which
are meromorphic along $D$, and it doesn't depend on $\alpha$. Note that the $F_{\alpha}$ may all be considered
as subsheaves of $F_{\infty}$.

We can define a tensor product of torsion-free parabolic sheaves: set
$$
(F\otimes G)_{\alpha}
$$
to be the subsheaf of $F_{\infty}\otimes _{\cO _X}G_{\infty}$ generated by the $F_{\beta '} \otimes G_{\beta
''}$ for $\beta '+ \beta '' \leq \alpha$.

On the other hand, if $E$ is a torsion-free sheaf on $X$ then it may be considered as a parabolic sheaf
(we say ``with trivial parabolic structure'') by setting $E_{\alpha}$ to be $E(\sum a_iD_i)$ for
$a_i$ the greatest integer $\leq \alpha _i$.

With this notation we may define for any vector bundle $E$ on $X$ the parabolic bundle
$$
E(\sum _{i=1}^k \alpha _i D_i) := E\otimes \cO _X(\sum _{i=1}^k \alpha _i D).
$$

\begin{lemma}
\label{structure}
Any parabolic line bundle has the form $L(\sum _{i=1}^k \alpha _i D_i)$ for $L$ a line bundle on $X$.
This may be viewed as $L(B)$ where $B$ is a rational divisor on $X$ (supported on $D$).
\end{lemma}
\eop

\begin{definition}
\label{lapb}
A parabolic sheaf $F$ is a {\em locally abelian parabolic bundle} if, in a Zariski neighborhood of any point
$x\in X$ there is an isomorphism between $F$ and a direct sum of parabolic line bundles.
\end{definition}

Most questions about locally abelian parabolic bundles can be treated by reducing to a local question then
looking at the case of
line bundles and using the structure result of Lemma \ref{structure}.

\subsection{Relationship with bundles on DM-stacks}

This material is the subject of a number of recent papers and preprints by Biswas \cite{Biswas2},
Cadman \cite{Cadman}, Matsuki and Olsson \cite{MatsukiOlsson} and specially Niels Borne \cite{Borne}.
\footnote{We thank A. Chiodo for useful discussions about this and for informing us of these references.}  

Suppose $f:(X',D')\rightarrow (X,D)$ is a morphism of smooth varieties with normal crossings divisors, such
that $f^{-1}(D)\subset D'$. Then we would like to define the {\em pullback} $f^{\ast}F$ of a locally abelian
parabolic bundle
$F$ on $(X,D)$, as a locally abelian parabolic bundle on $(X',D')$. This is not entirely trivial to do.
We propose two approaches. We first will pass through the notion of bundles on Deligne-Mumford stacks in order to
give a precise definition. On the other hand, we can use the locally abelian structure to reduce to the
case of parabolic line bundles. In this case we require that
$$
f^{\ast}(L(B)) = (f^{\ast}L)(f^{\ast}B)
$$
using the pullback $f^{\ast}(B)$ for a rational divisor $B$ (note that if $B$ is supported on $D$ then
the pullback will be supported on $D'$).

Any construction using a local choice of frame as in Definition \ref{lapb} or something else such as the
resolution we will discuss later, leaves open the problem of seeing that the construction is independent
of the choices which were made. In order to get around this kind of problem, we use the relationship
between parabolic bundles and bundles on certain Deligne-Mumford stacks. 

We define a DM-stack denoted $Z:= X[ \frac{D_1}{n_1}, \ldots , \frac{D_k}{n_k}]$.
Localizing in the etale topology over $X$
we may assume that the $D_i$ are defined by equations $z_i$. Only some of the components will
appear in any local chart; renumber things so that these are $z_i=0$ for $i=1,\ldots , k'$.
Then define the local chart for $Z$ to be given with coordinates $u_i$ by the equations $z_i = u_i^{n_i}$
for $i=1,\ldots , k'$ and $z_i=u_i$ for the other $i$. Without repeating the general theory of DM-stacks,
this defines a smooth DM-stack $Z$ whose coarse moduli space is $X$ and whose
stabilizer groups along the components $D_i$ are the groups $\Z / n_i \Z$.

This DM-stack may alternatively be considered as a ``Q-variety'' in the sense of Mumford \cite[section 2]{Mu}.

On $Z$ we have divisors which we may write as $\frac{D_i}{n_i}$, given in the local coordinate patch by
$u_i=0$. In particular, if $B$ is a rational divisor
supported along $D$ such that $B=\sum b_iD_i$ and if $n_ib_i\in \Z$ (that is, the denominator of $b_i$ divides
$n_i$) then $B$ becomes an actual divisor on $Z$ by writing
$$
B= \sum (n_ib_i) \frac{D_i}{n_i}.
$$
Let $p:Z\rightarrow X$ denote the projection.
If $E$ is a vector bundle on $Z$ then $p_{\ast}E$ is a torsion-free sheaf on $X$, and in fact it is a bundle.
We may define the {\em associated parabolic bundle ${\bf a}(E)$ on $(X,D)$} by the formula
$$
{\bf a}(E)_{\alpha} := p_{\ast} (E (\sum _i \alpha '_iD_i))
$$
where $\alpha '_i$ is the greatest rational number $\leq \alpha _i$ with denominator dividing $n_i$.

\begin{lemma}
\label{equivParaDM}
The above construction establishes an equivalence of categories between vector bundles on $Z$,
and locally abelian parabolic bundles on $(X,D)$ whose weights have denominators dividing the $n_i$.
\end{lemma}
\begin{proof}
See \cite[Theorem 5]{Borne}.  The setup there is slightly different in that all of the components 
$D_i$ are combined together, and Borne considers more generally torsion-free coherent sheaves. 
To adapt this to our situation, 
we will just comment on why the parabolic bundle on $X$ obtained from a bundle on $Z$
must satisfy the
locally abelian condition. Suppose $E$ is a bundle on the DM stack $Z$. It suffices to see the following claim: 
that there is a Zariski
open covering $X=\bigcup _i X_i$ which induces a covering $Z=\bigcup _i Z_i$ with $Z_i := Z\times _XX_i$,
such that every $E|_{Z_i}$ splits as a direct sum of line bundles on the DM-stacks $Z_i$. 
Fix a point $x\in X$ on a crossing point of components which we may assume are numbered $D_1,\ldots , D_b$. 
In a Zariski neighborhood $X_i$ of $x$, we have by construction
$$
Z_i = Y_i / G_i
$$
where $Y_i$ is a smooth scheme, $G_i = \Z / n_1 \times \cdots \times \Z / n_b$, 
$$
Y_i \rightarrow X_i
$$
is a Galois ramified covering with group $G_i$ and $X_i$ is the categorical quotient,
and $y\in Y_i$ is the unique point lying over $x$ and $G_i \cdot y = y$. 
The bundle $E|_{Z_i}$ is given by the data of a $G_i$-equivariant bundle $E_{Y_i}$ on $Y_i$.
In particular $G_i$ acts on the fiber $E_{Y_i, y}$ and since $G_i$ is abelian, there is a direct sum decomposition
$$
E_{Y_i, y} = \bigoplus _{j= 1}^r \comx (\xi _j)
$$
where $\xi _j$ are characters of $G_i$. Define the $G_i$-equivariant direct sum of line bundles
$$
F := \bigoplus _{j= 1}^r \cO _{Y_i}(\xi _j) := \bigoplus _{j= 1}^r \cO _{Y_i}\otimes _{\comx} \comx (\xi _j).
$$
Let $u_0: E_{Y_i}\rightarrow F$ be any morphism which restricts over $y$ to the given 
$G_i$-equivariant isomorphism
$$
E_{Y_i, y}\cong F_y,
$$
and put $u := \sum _{g\in G_i}g^{-1}u_0\circ g$. This is now $G_i$-equivariant, and is still an isomorphism over $y$.
Let $Y'_i$ be the Zariski open set where $u$ is an isomorphism. It is $G_i$-equivariant so it comes from an open
neighborhood $X_i'$ of $x$. Over the corresponding stack $Z'_i$ our morphism $u$ descends to an isomorphism between $E$
and a direct sum of line bundles. Thus the parabolic bundle on $X$ will be a direct sum of parabolic line bundles 
locally in the Zariski topology. We refer to \cite{Borne} for the remainder of the proof. 
\end{proof}

One expects more general parabolic bundles to correspond
to certain saturated sheaves on these DM stacks but that goes beyond our present requirements.

The inverse functor in Lemma \ref{equivParaDM} will be taken as the definition of the pullback $p^{\ast}$ from certain
locally abelian parabolic bundles on $(X,D)$ (that is, those with appropriate denominators)
to bundles considered as having their trivial parabolic structure on the stack $Z$ via the map $p:Z\rightarrow X$.

\begin{lemma}
\label{indep1}
The inverse functor $p^{\ast}$ is compatible with morphisms of DM-stacks of the form
$$
Z_1:=X[ \frac{D_1}{n_1}, \ldots , \frac{D_k}{n_k}]
\stackrel{q}{\rightarrow }
Z_2:=
X[ \frac{D_1}{m_1}, \ldots , \frac{D_k}{m_k}]
$$
whenever $m_i$ divides $n_i$. In other words, if $F$ is a parabolic bundle whose denominators divide the $m_i$
and if $E=p^{\ast}F$ is the associated bundle on $Z_2$, then $q^{\ast}E$ is the associated bundle on $Z_1$
(by a natural isomorphism satisfying a cocycle condition).
\end{lemma}
\eop

\begin{lemma}\label{pullback}
Suppose $f:(X',D')\rightarrow (X,D)$ is a morphism of smooth varieties with normal crossings divisors, such
that $f^{-1}(D)\subset D'$. Then the {\em pullback} $f^{\ast}F$ of a locally abelian
parabolic bundle
$F$ on $(X,D)$, is defined as a locally abelian parabolic bundle on $(X',D')$.
\end{lemma}
\begin{proof}
Using the equivalence of Lemma \ref{equivParaDM}, we can define the pullback of a locally abelian parabolic bundle.
Indeed, given a morphism
$f:(X',D')\rightarrow (X,D)$ and any positive integers $n_1,\ldots , n_k$ for the components $D_i$ of $D$,
there exist integers $n'_1,\ldots ,n'_{k'}$ for the components $D'_i$ of $D'$ such that $f$ extends to
a morphism
$$
f_{DM}:X' [ \frac{D'_1}{n'_1}, \ldots , \frac{D'_{k'}}{n'_{k'}}]
\rightarrow
X[ \frac{D_1}{n_1}, \ldots , \frac{D_k}{n_k}]
$$
If $F$ is a locally abelian parabolic bundle on $X$ we can choose integers $n_i$ divisible by the 
denominators of the weights of $F$,
then choose $n'_j$ as above. Thus $F$ corresponds to a bundle $E$ on the DM-stack
$X[ \frac{D_1}{n_1}, \ldots , \frac{D_k}{n_k}]$. The pullback $f^{\ast}_{DM}(E)$ is a bundle on
the DM-stack $X'[ \frac{D'_1}{n'_1}, \ldots , \frac{D'_{k'}}{n'_{k'}}]$.
Define $f^{\ast}(F)$ to be the locally abelian parabolic bundle on $(X',D')$ corresponding
(again by Lemma \ref{equivParaDM}) to
the bundle $f^{\ast}_{DM}(E)$.

Using Lemma \ref{indep1}, the parabolic bundle $f^{\ast}(F)$ is independent of the choice of $n_i$ and $n'_j$.
\end{proof}

This definition can be extended, using local charts which are schemes, to the pullback to any DM-stack,
and it is compatible with the previous notation in the sense that if $p:Z\rightarrow X$ is
the projection from the DM-stack used in Lemma \ref{equivParaDM} to the original variety, then
the bundle $E=p^{\ast}(F)$ which corresponds to $F$ is indeed the pullback as defined two paragraphs ago
(which would be a parabolic bundle on the DM-stack $Z$ but which is actually a regular bundle with trivial parabolic
structure).

\subsection{Chern character of parabolic bundles}\label{ChernPar}

Mumford, Gillet, Vistoli (\cite{Mu},\cite{Gi},\cite{Vi}) have defined Chow groups for DM-stacks. Starting with $(X,D)$ 
and choosing denominators $n_i$ we obtain the
stack $Z:= X[ \frac{D_1}{n_1}, \ldots , \frac{D_k}{n_k}]$. The coarse moduli space of $Z$ is the original $X$,
so from \cite[Theorem 6.8]{Gi}, the pullback and pushforward maps establish an isomorphism of rational Chow groups
$$
CH^{\ast}(Z)_{\Q }\cong CH^{\ast}(X)_{\Q }.
$$
We also have a Chern character for bundles on a DM-stack \cite[section 8]{Gi}.
If $F$ is a locally abelian parabolic bundle on $(X,D)$, choose $n_i$ divisible by
the denominators of the weights of $F$ and define the Chern character of $F$ by
$$
{\rm ch}^{\rm par}(F):= p_{\ast} ({\rm ch}(p^{\ast}F)) \in CH^{\ast}(X)_{\Q }.
$$
This doesn't depend on the choice of $n_i$, by Lemma \ref{indep1}, and the compatibility of the Chern character with
pullbacks on DM-stacks.

Using the fact that Chern character commutes with pullback for bundles on DM-stacks, we obtain:

\begin{lemma}
\label{pullbackCh}
Suppose $f:(X',D')\rightarrow (X,D)$ is a morphism of smooth varieties with normal-crossings divisors,
such that $f^{-1}(D)$ is supported on $D'$. If $F$ is a locally abelian parabolic bundle on $(X,D)$ then
the pullback $f^{\ast}F$ is a locally abelian parabolic bundle on $(X',D')$ and this construction commutes with
Chern character:
$$
{\rm ch}^{\rm par}(f^{\ast}F) = f^{\ast}({\rm ch}^{\rm par}(F))\in CH^{\ast}(X')_{\Q }.
$$
\end{lemma}
\eop

If $F$ is a usual vector bundle considered with its trivial parabolic structure, then ${\rm ch}(F)$ is the usual
Chern character of $F$ (this is the case $n_i=1$ in the definition). Also, if $f$ is a morphism then the
pullback $f^{\ast}(F)$ is again a usual vector bundle with trivial parabolic structure.

We can easily describe the pullback and Chern characters for parabolic line bundles. If $B = \sum b_iD_i$ is
a rational divisor on $X$, supported on $D$, and if $f:(X',D')\rightarrow (X,D)$ is a morphism as before,
we have the formula
$$
f^{\ast}(\cO _X(B)) = \cO _{X'}(B')
$$
where $B'=\sum b_i (f^{\ast} D_i)$ with $f^{\ast}(D_i)$ the usual pullback of Cartier divisors.
Similarly, if $E$ is a vector bundle on $X$ then we have
$$
f^{\ast}(E (B)) = (f^{\ast}E)(B').
$$

For the Chern character, recall that ${\rm ch}(\cO _X(B))=e^B$ for a divisor with integer coefficients $B$.
This formula extends to the case of a rational divisor, to give the formula for the Chern
character of a parabolic line bundle in the rational Chow group.
More generally for a vector bundle twisted by a rational divisor we have
$$
{\rm ch} ^{\rm par}(E(B)) = {\rm ch}(E)\cdot e^B.
$$

The Chern character is additive on the $K_0$-group of vector bundles on a DM-stack, so (modulo saying something
good about choosing appropriate denominators whenever we want to apply the formula) we get the same
statement for the Chern character on the $K_0$-group of locally abelian parabolic bundles. One should say what is
an {\em exact sequence} of locally abelian parabolic bundles. Of course a morphism of parabolic sheaves $F\rightarrow G$
is just a collection of morphisms of sheaves $F_{\alpha}\rightarrow G_{\alpha}$ and the kernel and cokernel are still
parabolic sheaves.  
A short sequence
$$
0\rightarrow F \rightarrow G \rightarrow H\rightarrow 0
$$
of parabolic sheaves is {\em exact} if the resulting sequences
$$
0\rightarrow F_{\alpha} \rightarrow G_{\alpha} \rightarrow H_{\alpha}\rightarrow 0
$$
are exact for all $\alpha$. This notion is preserved by the equivalence of Lemma \ref{equivParaDM} and indeed
more generally the pullback functor $f^{\ast}$ preserves exactness.
All objects in these sequences are bundles so there are no higher $Ext$ sheaves and in particular there
is no requirement of flatness when we say that  $f^{\ast}$ preserves exactness, although the $p^{\ast}$ of
Lemma  \ref{equivParaDM} is flat.

Thus we can define the $K_0$ of
the category of locally abelian parabolic bundles, and the Chern character is additive here.
Also the Chern character is multiplicative on tensor products (this might need some further proof which we don't
supply here). And the pullback $f^{\ast}$ is defined on the $K_0$ group.

{\em Caution:} The Chern character doesn't provide an isomorphism between $K_0$ and the Chow group, on a DM-stack.
For this one must consider an extended Chern character containing further information over the stacky locus,
see \cite{Toen} \cite{ChenRuan} \cite{AbramovichGraberVistoli}. Thus, the same holds for parabolic bundles:
two parabolic bundles with the same Chern character are not necessarily the same in the $K_0$-group of parabolic
bundles.

We don't give a formula for the Chern character of a parabolic bundle $F$ in terms of the
Chern characters of the constituent bundles $F_{\alpha}$.\footnote{This point seems to be somewhat delicate, for example
D. Panov in his thesis \cite{Panov} gives corrected versions of some of the formulae obtained by Biswas \cite{Biswas}.}
The formula we had claimed in the first version of the present paper was wrong. We will treat this question elsewhere.
For our present purposes, what we need to know is contained in the following lemma and its corollary.

\begin{lemma}
\label{diffpar}
Suppose $F$ is a parabolic bundle. Then the difference between $F$ and any one of its $F_{\alpha}$ in the rational Chow ring
of $X$ or equivalently the rational $K$-theory of $X$,
is an element which is concentrated over $D$ (i.e. it may be represented by a rational combination of sheaves concentrated on $D$).
\end{lemma}
\begin{proof}
There is a localization sequence in rational $K$-theory of $X$ (see \cite[Proposition 5.15]{Sr}) :
$$
K_0(D)\otimes \Q \rightarrow K_0(X)\otimes \Q \rightarrow K_0(X-D)\otimes\Q \rightarrow 0.
$$
The restrictions of $F$ and $F_{\alpha}$ to $U=X-D$ are isomorphic vector bundles so ${\rm ch}(F)-{\rm ch}(F_{\alpha})$
maps to zero in $K_0(X-D)\otimes\Q $. Thus it comes from $K_0(D)\otimes \Q$.
\end{proof}

\begin{corollary}\label{co.-ch}
If $ch^{\rm par}(F)\in CH^0(X)_{\Q}$ for a parabolic bundle, 
then any $F_{\alpha}$ is equivalent to an element coming from $CH^{\ast}(D)$.
\end{corollary}

\section{The parabolic bundle associated to a logarithmic connection}\label{Parlog}

Suppose as above that $X$ is a smooth projective variety with normal crossings
divisor $D$, and  let $U:= X-D$. Write as before
$D=D_1\cup \cdots \cup D_k$ the decomposition of $D$ into smooth irreducible
components.

Suppose $(E,\nabla )$ is a vector bundle with logarithmic connection $\nabla$ on $X$
such that
the singularities of $\nabla$ are concentrated over $D$. Fix a component $D_i$ of the divisor.
Define the {\em residue} of $(E,\nabla )$ along $D_i$ to be the pair $(E|_{D_i}, \eta _i)$
where  $E|_{D_i}$ is the restriction of $E$ to $D_i$.
The residual transformation $\eta _i$ is
the action of the vector field $z\frac{\partial}{\partial z}$ where $z$ is a local coordinate of $D_i$.
This is independent of the choice of local coordinate. 
The connection induces an operator
$$
\nabla _{D_i}: E|_{D_i}\rightarrow (E|_{D_i})\otimes _{\cO _{D_i}} (\Omega ^1_X(\log D))|_{D_i},
$$
whose projection by the residue map is $\eta _i$. 
This will not in general lift to a connection on $E|_{D_i}$ (an error in our first version pointed out by 
H. Esnault). Such a lift exists locally in the etale topology where we can write $X$ as a product of $D_i$ and 
the affine line, and in this case $\eta _i$ is an automorphism of the  bundle with
connection $E|_{D_i}$, in particular the eigenvalues of $\eta _i$ are locally constant functions
along $D_i$. The latter holds without etale localization, and the eigenvalues are
actually constant because we assume $D_i$ irreducible.

Assume the following condition:

\begin{definition}
\label{rationalresidues}
We say that $(E,\nabla )$ {\em has rational residues} if
the eigenvalues of the residual transformations of $\nabla$ along components of $D$ (that is,
the $\eta _i$ above) are rational numbers.
\end{definition}

In this case we will construct a parabolic vector bundle $F=\{ F_{\alpha }\}$ on $(X,D)$ together with
isomorphisms
$$
F_{\alpha}|_{U} \cong E|_{U}.
$$
We require that $\nabla$ extend to a logarithmic connection on each $F_{\alpha}$. Finally, we require that
the residue of $\nabla$ on the piece $F_{\alpha} / F_{\alpha - \varepsilon \delta ^i}$ which is concentrated over $D_i$,
be an operator with eigenvalue $-\alpha_i$.  

\begin{definition}
\label{apb}
A parabolic bundle with these data and properties is
called a {\em parabolic bundle associated to $(E,\nabla )$}.
\end{definition}

The notion of a parabolic bundle associated to a connection was discussed in \cite{Inaba} and 
\cite{InabaEtAl}. In those places, no restriction is placed on the residues, and parabolic structures with arbitrary
complex numbers for weights are considered. They use full flags and treat the case of curves.
It isn't immediately clear how the full flags would generalize for normal crossings divisors in higher dimensions.

\begin{lemma}\label{le.-parabolic}
For any vector bundle with logarithmic connection having rational residues, 
there exists a unique associated parabolic bundle, and furthermore it is locally abelian.
\end{lemma}
\begin{proof} 
This construction is basically the one discussed by Deligne in \cite{De} Proposition 5.4 and attributed to
Manin \cite{Manin}. We give three constructions. 

\noindent
{\bf (I)}\, 
Suppose $E$ is a vector bundle with a logarithmic
connection on $(X,D)$. Applying the discussion of \cite{De} we may assume the eigenvalues are
contained in $(-1,0]$.  Consider the residue $\eta_i \in
\mbox{End}(E_{|D_i})$ whose eigenvalues $-\al^j_i$ are rational
numbers and satisfying $\al^1_i<\al^2_i<...<\al^{n_i}_i$.

The eigenspaces of the fibres of $E_{|D_i}$ given by the endomorphism $\eta_i$ associated to an eigenvalue 
$-\al^j_i$ (which is a constant), defines
a subsheaf $A_{\al^j_i}$ of $E_{|D_i}$, called the \textit{generalized eigenspaces}. 
Let
$F^{\al^j_i}_i= \sum_{l=1}^jA_{\al^l_i}$. Define the subsheaves
$\ov{F^{\al^j_i}}$ of $E$ by the exact sequence:
$$0\lrar \ov{F^{\al^j_i}_i}\lrar E\lrar E_{|D_i}/F^{\al^j_i}_i\lrar 0. $$

Given a set of rational weights $\al=(\al_1,\al_2,...,\al_k)$, i.e, $\al_i$
is a weight along $D_i$, we can associate a subsheaf
$$
F_{\al}=\bigcap _{i=1}^k\ov{F^{\al_i}_i}
$$
of $E$. This defines a parabolic structure on $(X,D)$. 
Moreover, each
$F_\al$ restricts on $U$ to $E_{|U}$.

Since the residue endomorphism $\eta_i$ preserves the generalized eigenspaces, it preserves the sheaves $F^{\al^j_i}_i$ and
hence $\nabla$ induces a logarithmic connection on
each $\ov{F^{\alpha^j_i}_i}$ and thus on $F_\al$.

For the locally abelian condition, note that
the same construction may be done over the analytic topology, and since the local monodromy groups are abelian
on punctured neighborhoods of the crossing points of $D$ which are products of punctured discs, we get a 
decomposition into a direct sum of parabolic line bundles, locally in the analytic topology. 
Artin approximation gives it locally in the etale topology,
and the argument mentioned in Lemma \ref{equivParaDM} gives the locally abelian condition in the Zariski topology.

\noindent
{\bf (II)} \,
In \cite[Proposition 5.4]{De} Deligne constructs a bundle associated to any choice of lifting function 
$\tau : \comx / \Z \rightarrow \comx $. In fact the same construction will work if we choose different liftings $\tau _i$ for
each component of the divisor $D_i$. Given $\alpha$, choose liftings $\tau _i$ which send $\R / \Z$ to the intervals
$[-\alpha _i,1 -\alpha _i )$ respectively. Applying the extension construction gives the bundle $F_{\alpha}$ and
these organize into our parabolic bundle on $X$, locally abelian by the same reasoning as in the previous paragraph.

\noindent
{\bf (III)} \,
Choose $n_i$ such that the residues of
$\nabla $ along $D_i$ are integer multiples of $\frac{1}{n_i}$. Then the parabolic bundle associated to 
$(E,\nabla )$ corresponds to the Deligne canonical extension of $(E,\nabla )$ over the Deligne-Mumford stack
$Z:= X[\frac{D_1}{n_1},\ldots , \frac{D_k}{n_k}]$, see \S 2.2. Note that the pullback of $(E,\nabla )$ to a logarithmic
connection over this Deligne-Mumford stack, has integer residues so exactly the construction of \cite{De}
can be carried out here. Uniqueness of the construction implies descent from
an etale covering of the stack to the stack $Z$ itself. As noted in the proof of Lemma \ref{equivParaDM},
a bundle on $Z$ corresponds to a locally abelian parabolic structure on $X$. The fact that the bundle has a
logarithmic connection over $Z$ translates into the statement that the $F_{\alpha}$ are preserved by the logarithmic
connection which
is generically defined from $U$. 
\end{proof}

\begin{lemma}
\label{eigenvalues}
Let $F$ be the parabolic bundle associated to a logarithmic connection. 
On the bundle $F_{\alpha}$ the eigenvalues of the residue of $\nabla$ along $D_i$ are
contained in the interval $[- \alpha _i, 1-\alpha _i )$.
\end{lemma}
\begin{proof}
The restriction to $D_i$ may be expressed as 
$$
F_{\alpha}|_{D_i} = F_{\alpha} / F_{\alpha - \delta ^i},
$$ 
so the condition in the definition of associated parabolic bundle fixes the eigenvalues on the graded pieces 
here as being in the required interval. 
\end{proof}

The proof of \ref{le.-parabolic} also shows that
\begin{lemma}\label{le.-parabolic2}
Given an exact sequence of logarithmic connections with rational
residues along $D$ :
$$0\lrar E''\lrar E\lrar E'\lrar 0$$
there is an exact sequence of parabolic bundles
$$0\lrar F_{E''}\lrar F_{E}\lrar F_{E'}\lrar 0$$
such that $F_{E''}, F_{E}, F_{E'}$ are the parabolic bundles
associated to $E'',E$ and $E'$ respectively. Similarly the construction is compatible with direct sums. 
\end{lemma}

\begin{lemma}
\label{commuterestrict}
Suppose $Y\subset X$ is a subvariety intersecting $D$ transversally. Let $F$ be the parabolic bundle associated to a
vector bundle with logarithmic connection $(E,\nabla )$ over $X$. Then $F|_Y$ is the parabolic bundle associated to
the restriction of $(E,\nabla )$ to $Y$.
\end{lemma}
\begin{proof}
Since the intersection of $Y$ and $D$ is transversal, $D$ restricts to a normal crossing divisor $D_Y$ on $Y$. 
The restriction of $(E,\nabla)$ on $Y$ has rational residues $\{\al_i\}$ which is a subset of the rational residues of
$(E,\nabla)$ and corresponds to those sheaves $F_{\al_i}$ which restricted to $Y$ are nonzero. 
In other words, using the construction of parabolic structure
in Lemma \ref{le.-parabolic}, it follows that $F{|}_Y$ is the parabolic bundle
associated to the restriction of $(E,\nabla)$ to $Y$.
\end{proof}

\subsection{Computing cohomology using the associated parabolic bundle}

Suppose $X$ is a smooth projective variety  with normal crossings divisor $D\subset X$ and let $U:= X-D$.
Let $j: U\rightarrow X$ denote the inclusion.

Suppose $(E,\nabla )$ is a vector bundle on $X$
with logarithmic connection and rational reesidues along $D$, with associated parabolic bundle $F$.
Let $DR(E ,\nabla )$ denote the de Rham complex
$$
E \rightarrow E \otimes _{\cO _X}\Omega ^1_X(\log D) \rightarrow
E \otimes _{\cO _X}\Omega ^2_X(\log D) \rightarrow \cdots
$$
with $d_{\nabla}$ as differential. 

Define the {\em logarithmic de Rham cohomology} 
$$
\cH ^{\cdot} _{DR}(X; E,\nabla ) := \bbH ^{\cdot}(X, DR(E ,\nabla )).
$$
In this notation the fact that it is logarithmic along $D$ is encoded in the fact that $\nabla$ is
a logarithmic connection along $D$. On the other hand, let $E_U^{\nabla}$ denote the local system
on $U$ of flat sections for the connection $\nabla$. 

\begin{theorem}\label{logcoh}
{\rm (Deligne)}
Suppose that the residues of $\nabla$ have no strictly positive integers as eigenvalues. Then
$$
\cH ^{i} _{DR}(X; E,\nabla )\cong H^i(U, E_U^{\nabla}).
$$
\end{theorem}
\begin{proof}
This is \cite{De}, Proposition 3.13 and Corollaire 3.14.
\end{proof}

\begin{corollary}
\label{logcohpar}
Suppose $(E,\nabla )$ is a logarithmic connection. Let $F$ denote the associated parabolic bundle,
so that for any $\alpha$ we obtain again a logarithmic connection $(F_{\alpha}, \nabla )$. 
If $\alpha = (\alpha _1,\ldots , \alpha _k)$
with all $\alpha _i >0$ then 
$$
\cH ^{i} _{DR}(X; F_\alpha,\nabla ) \cong H^i(U, E_U^{\nabla}).
$$
\end{corollary}
\begin{proof}
By Lemma \ref{eigenvalues} the eigenvalues of the residue of $\nabla$ on $F_{\alpha}$ are in
$[- \alpha _i, 1-\alpha _i )$. If $\alpha _i > 0$ this is guaranteed not to contain any positive integers. 
\end{proof}

\subsection{The relative logarithmic de Rham complex}

Suppose $(E,\nabla )$ is a vector bundle with logarithmic connection relative to a pair $(X,D)$ where $X$ is smooth and
$D$ is a normal crossings divisor. In practice, $E$ will be of the form $F_{\alpha}$ as per Theorem \ref{logcohpar}.

Suppose we have a map $f:X\rightarrow Y$ such that components of $D$ are either
over the singular divisor $J$ in $Y$, or have relative normal crossings over $Y$ \cite{De}.
We have the logarithmic complex $\Omega ^{\cdot}_X(\log D)$ and the subcomplex
$f^{\ast} \Omega ^{\cdot}_Y(\log J)$. Put
$$
\Omega ^{\cdot}_{X/Y}(\log D):= \Omega ^{\cdot}_X(\log D)/f^{\ast} \Omega ^{\cdot}_Y(\log J).
$$
We can form a differential $d_{\nabla}$ using $\nabla$ on
$E\otimes \Omega ^{\cdot}_{X/Y}(\log D)$. Call the resulting logarithmic de Rham complex $DR(X/Y, E,\nabla )$.  In this
notation the divisors $D$ and $J$ are implicit. 

Define a complex of quasicoherent sheaves on $Y$ as follows:
$$
\cH ^{\cdot} _{DR}(X/Y; E,\nabla ) := \R ^{\cdot} f_{\ast}DR(X/Y, E,\nabla ).
$$
To define this precisely, choose an open covering of $X$ by affine open sets $X_{\beta}$ such that the map $f$
is affine on the multiple intersections. Let $\check{\cC}^{\cdot}DR(X/Y, E,\nabla )$ denote the 
simple complex associated to the double complex obtained by applying the \v{C}ech complex in each degree.
The components are acyclic for $f_{\ast}$ so we can set
$$
\R ^{\cdot} f_{\ast}DR(X/Y, E,\nabla ):= f_{\ast}\check{\cC}^{\cdot}DR(X/Y, E,\nabla ).
$$
The component sheaves of the \v{C}ech complex on $X$ are quasicoherent but the differentials are not $\cO_X$-linear.
However, the differentials are $f^{-1}\cO _Y$-linear, so the direct image complex consists of quasicoherent sheaves on
$Y$. 

The complex $DR(X/Y, E,\nabla )$ is filtered by the Hodge filtration which is the ``filtration b\^ete''.
This presents $DR(X/Y, E,\nabla )$ as a successive extension of complexes $E\otimes \Omega ^{m}_{X/Y}(\log D)[-m]$.
The higher derived direct image becomes a successive extension of the higher derived direct images of these
component complexes, for example this is exactly true if we use the \v{C}ech construction above.

The standard argument of Mumford \cite{Mu} shows that each
$$
\R ^{\cdot} f_{\ast}E\otimes \Omega ^{m}_{X/Y}(\log D)[-m]
$$
is a perfect complex over $Y$ and formation of this commutes with base changes $b:Y'\rightarrow Y$. 
Therefore the same is true of the successive extensions. Thus formation of
$\cH ^{\cdot} _{DR}(X/Y; E,\nabla )$ commutes with base change. 

At this point we are imprecise about the conditions on $Y'$ and $b$. All objects should be defined after base change.
For example if $X':= Y'\times_Y X$ is smooth and the pullbacks of $D$ to $X'$ and $J$ to $Y'$ are normal crossings,
then everything is still defined. In \S 6 we will want to consider a more general situation in which the logarithmic
de Rham complex is still defined. To be completely general, refer to Illusie \cite{Illusie} for
example for an intrinsic definition, but we could take as definition
$$
DR(X'/Y', E',\nabla '):= pr_2^{\ast} DR(X/Y, E,\nabla )
$$
for $pr_2: X'\rightarrow X$ the projection, and with that the base-change formula holds for any $b:Y'\rightarrow Y$. 

This gives the following lemma which is well known.

\begin{lemma}
\label{reldrperfect}
Suppose $(E,\nabla )$ is a vector bundle with logarithmic connection on $(X,D)$. 
The relative logarithmic de Rham cohomology  $\cH ^{\cdot} _{DR}(X/Y; E,\nabla )$ is a perfect complex on $Y$
and if $b:Y'\rightarrow Y$ is a morphism then letting $X':= Y'\times_Y X$ and $(E',\nabla ')$ be the pullbacks, we have
a quasiisomorphism
$$
\cH ^{\cdot} _{DR}(X'/Y'; E',\nabla ') \stackrel{{\rm qi}}{\sim} b^{\ast} \cH ^{\cdot} _{DR}(X/Y; E,\nabla ).
$$
\end{lemma}

The complex $\cH ^{\cdot} _{DR}(X/Y;  E,\nabla )$ has a logarithmic 
Gauss-Manin connection which induces the usual Gauss-Manin connection
on the cohomology sheaves. This is well known from the theory of $\cD$-modules, see also \cite{KatzICM}. We discuss it in \S 6. 

As Deligne says in \cite{De} 3.16, the calculation of cohomology also works in the relative case. 
Restrict over the 
open set $U$ which is the complement of $J$. Then $X_U\rightarrow U$ is a smooth map and the divisor $D_U:= D\cap X_U$ is a union of
components all of whose multiple intersections are
smooth over $Y$.  Let $W:= X_U-D_U$. Let $f_W:W\rightarrow U$ be the map restricted to our open set. With the current hypotheses,
topologically it is a fibration. The relative version
of Lemma \ref{logcoh} (see also \cite{KatzICM}) says the following. 

\begin{lemma}
\label{logcohrel}
In the above situation, suppose that the eigenvalues of the residue of $\nabla$ along horizontal
components, that is components of $D_U$, are never positive integers.  Then 
the vector bundle with Gauss-Manin connection $\cH ^{i} _{DR}(X_U/U; E_U,\nabla )$ on
$U$ is the unique bundle with connection and regular singularities corresponding to the local system
$R^i f_{W,\ast} \comx _W$ on $U$. 
\end{lemma}

\subsection{Chern characters for higher direct images}

We now introduce some notation for calculations in Chow groups.  Suppose $f: X\rightarrow Y$ is a flat morphism.
If $E$ is a vector bundle on $X$
then its higher direct image complex $\R ^{\cdot} f_{\ast}(E)$ is a perfect complex on $Y$ and we can define its Chern character
by
$$
{\rm ch}(\R ^{\cdot} f_{\ast}(E)) := \sum _i (-1)^i  {\rm ch}(\R ^{i} f_{\ast}(E)).
$$
The Grothendieck-Riemann-Roch theorem implies the following observation which says that the Chern character of the  
higher direct image
complex is a function of the Chern character of $E$. 

\begin{proposition}\label{pr.-GRR}
Suppose $f:X\rightarrow Y$ is a proper morphism between smooth varieties. Then there is a map on rational Chow groups
$$
\chi  _{X/Y} : CH^{\cdot} (X)_\Q \rightarrow CH^{\cdot} (Y)_\Q
$$
which represents the Euler characteristic for higher direct images in the sense that
$$
{\rm ch}(\R ^{\cdot} f_{\ast}(E)) = \chi  _{X/Y} ({\rm ch}(E))
$$
for any vector bundle or coherent sheaf $E$ on $X$. 
\end{proposition}
\begin{proof}
The GRR formula gives an explicit formula for $\chi  _{X/Y}$ \cite[Theorem 15.2]{Fu}. 
One also notes that the Chow groups calculate the rational
$K_0$ and it is clear that ${\rm ch}\circ \R ^{\cdot} f_{\ast}$ is additive on exact sequences so it passes to a function 
on $K_0$.
\end{proof}

Esnault-Viehweg in \cite{Es-Vi2} use the explicit formula for $\chi  _{X/Y}$. We use less of GRR here in the sense that 
we use only existence of
the function rather than calculating with the formula.

\subsection{Modification of the base}

Consider the situation of a family $f:X\rightarrow S$. Let $\cO _X$ denote the trivial bundle with
its standard connection on $X$. 
We assume that $S$ is irreducible.
Let $U\subset S$ denote a nonempty open set over which $f$ is smooth, and let $X_U:= X\times _SU$.

The $\cH ^i_{DR}(X_U/U)$ are vector bundles on $U$. If the complementary divisor $J:= S-U$ has normal crossings then
the Gauss-Manin connections on the cohomology vector bundles have logarithmic singularities with rational residues (this
last statement is the content of the well-known ``monodromy theorem''). Let 
$$
\overline{\cH }^i_{DR}(X_U/U)
$$
denote the parabolic bundle on $S$ associated to this logarithmic connection. Our goal is to investigate the
alternating sum of the
Chern characters of these bundles on $S$. In view of this goal, we can make the following reduction.

\begin{lemma}
\label{modif-pullback}
In the above situation, suppose $p:S'\rightarrow S$ is a generically finite surjective map from an irreducible variety to $S$,
and if $U'\subset S'$ is a nonempty open set mapping to $U$, let $X':= X\times _SS'$. Assume that $J':= S'-U'$ again 
has normal crossings. Then the parabolic bundle $\overline{\cH }^i_{DR}(X_U/U)^{\rm par}$, which is defined as the
one which corresponds to the Gauss-Manin connection on $\cH ^i_{DR}(X'_{U'}/U')$, is isomorphic to the pullback
from $S$ to $S'$ of the parabolic bundle $\overline{\cH }^i_{DR}(X_U/U)^{\rm par}$ on $S$.
\end{lemma}
\begin{proof}
Notice that over $U'$ the statement is true for the cohomology bundle, since it is a base-change theorem for smooth morphisms. 
Now by Lemma \ref{pullback}, we deduce the statement for the parabolic bundle as well.
\end{proof}

\begin{corollary}
\label{modif}
Let $p:S'\rightarrow S$ and $U'$ be as in the previous lemma. Then 
$$
\sum _i (-1)^i {\rm ch}^{\rm par}(\overline{\cH }^i_{DR}(X_U/U))
= p^{\ast}
\sum _i (-1)^i {\rm ch}^{\rm par}(\overline{\cH }^i_{DR}(X_U/U)).
$$
In particular, $\sum _i (-1)^i {\rm ch}^{\rm par}(\overline{\cH }^i_{DR}(X_U/U))$ is in the degree zero
piece of the rational Chow group of $S$, if and only if 
$\sum _i (-1)^i {\rm ch}(\overline{\cH }^i_{DR}(X_U/U)^{\rm par})$ is in the degree zero
piece of the rational Chow group of $S'$.
\end{corollary}
\eop

Since our goal is to prove that the alternating sum in question at the end of the previous corollary is
in the degree zero piece of the rational Chow group, we can safely replace $X\rightarrow S$ and the open set $U$
by the base-change
$X'\rightarrow S'$ whenever $S'\rightarrow S$ is a generically finite rational map, with any nonempty open subset
$U'\subset S'$. In what follows we will often make this type of reduction. There are several different flavours,
for example we could simply decrease the size of the open set $U$; we could do any kind of birational transformation on $S$;
or we could take a finite covering of an open set of $S$ and complete to a smooth variety. 

In particular by using finite coverings we may assume that the monodromy transformations at infinity are unipotent
rather than just quasi-unipotent. In this case, the associated parabolic extension is just a regular bundle and it is
the same as the Deligne canonical extension. 

In order to avoid overly heavy notation in what follows, we will keep the same letters when we replace $S,X,U$
and speak of this as ``making a generically finite modification of $S$,'' or some similar such phrase.

\section{Lefschetz fibrations}

\subsection{Lefschetz pencil \cite{Ka}}

Suppose $X\subset \p^N$ is a nonsingular projective variety of dimension $n$, defined over the complex numbers.
Consider its dual variety $\hat{X}\subset \hat{\p^N}$. A Lefschetz pencil is a projective line 
$\p^1\subset \hat{\p^N}$ such that its intersection with $\hat{X}$
is a reduced zero-dimensional subscheme. Furthermore, these points correspond to hyperplane sections
of $X$ with only one ordinary double point as its singularity.

Then there is a Lefschetz fibration
\begin{eqnarray*}
Z\,\,\, &\sta{\alpha}{\lrar} & X \\
\downarrow {\scriptstyle \varphi} & & \\
\p^1.\,\, &&
\end{eqnarray*}
such that $Z$ is a blow--up of $X$ along the base locus $B$ defined by the vanishing of any two 
hyperplane sections $H_s,H_t\,\in \p^1$.

\begin{remark}\label{re.-mon}
We remark that if $n$ is odd, then the monodromy transformations of the
Lefschetz fibration in the middle dimensional cohomology, are reflections (i.e., they have order $2$ with 
one eigenvalue as $-1$ and the rest as $+1$). In paricular, they are not unipotent transformations.
\end{remark}

\subsection{Semistable reduction for families of rational curves}

We need the genus $0$ version of semistable reduction as proven by de Jong \cite{deJong}. 

\begin{lemma}
\label{sreduction}{\rm (de Jong)}
Suppose $P\rightarrow S$ is a morphism whose general fibers are projective lines, with a divisor $K\subset P$. 
Then there is a finite
covering and modification of the base $S'\rightarrow S$ and a family $P'\rightarrow S'$ with divisor $K'\subset P'$ such that
$(P',K')\rightarrow S'$ is a semistable family of marked rational curves, and such that the open set $U'\subset S'$ over which
the family is a smooth family of marked projective lines, has complement $D':= S'-U'$ which is a divisor with normal crossings.
Further, the variety $P'$ can be assumed to be smooth .
\end{lemma}
\begin{proof}
This is \cite[Theorem 4.1, Proposition 3.6]{deJong}.
\end{proof}

\subsection{A Lefschetz fibration for families}\label{sub.-lef}

Suppose $S$ is a smooth projective variety and $J\subset S$ is a normal crossings divisor with smooth components.
Let $U:= S-J$. Suppose $f:X_U\rightarrow U$ is a smooth projective morphism over $U$, where $X_U$ is smooth.
Let $n$ denote the relative dimension of $X_U$ over $U$. 

In the following discussion we allow modification of the base $S$ as described in  Lemma \ref{modif-pullback}
and Corollary \ref{modif}.

Suppose we are given an embedding $X_U\rightarrow S\times \p ^N$.
Choose a general line $\p^1\subset \hat{\p^N}$. After possibly restricting to a smaller open set we may assume that this
line defines a Lefschetz pencil for the subvariety $X_u$ for every $u\in U$.

Let $P_0:= U\times \p ^1$, 
with projection $$
q:P_0\lrar S
$$
so that we have a family of Lefschetz pencils
$$
X_U\sta{\alpha }{\longleftarrow} Z_0 \sta{\varphi}{\longrightarrow} P_0
$$
where
$\alpha : Z_0\longrightarrow X_U$ is blowing up the smooth family of base loci 
$B_0\longrightarrow U=S-J$, and the map $Z_0\longrightarrow P_0$ 
is a Lefschetz pencil of the fiber $X_s$
for $s\in S-J$.

Let $K_0\subset P_0$ denote the family of points over which the Lefschetz pencils have singular fibers. 
This is a divisor which is smooth over the base open set $U$: the different double points stay distinct as $u\in U$ moves
by the hypothesis that we have a Lefschetz pencil for every $u\in U$. 

By going to a finite covering of $U$ we may suppose that each irreducible component of the divisor $K_0\subset P_0$
maps isomorphically onto $U$. Thus we can consider $(P_0,K_0)$ as a family of marked rational curves.

Apply the semistable reduction result Lemma \ref{sreduction}, recalled in the previous subsection. 
After further modification of the base,
then completing the families $X_U$, $Z_0$ and $B_0$ we obtain the following situation. 

Keep the notations that 
$J \subset S$ is a divisor with normal crossings and $U:= S-J$ an open subset.
We have a semistable marked rational curve
$q: (P,K)\rightarrow S$ with divisor $D:= K\cup  q^{-1}(J)\subset P$. We have a diagram
$$
\begin{array}{ccc}
Z & \stackrel{\varphi}{\rightarrow} & P \\
\downarrow{{\scriptstyle \alpha}} & & \downarrow{{\scriptstyle q}}\\
X & \stackrel{f}{\rightarrow} & S.
\end{array}
$$
We also have a subvariety $B\subset X$ corresponding to the base loci.

If we denote in general by a subscript $U$ the restriction of objects over the open set $U$, then
$P_U\cong U\times \p ^1$ and 
 $$
X_U\sta{\alpha _U}{\longleftarrow} Z_U \sta{\varphi _U}{\longrightarrow} P_U\cong U\times \p ^1
$$
is a family of Lefschetz pencils indexed by $U$. The divisor $K_U\subset P_U$ which is a union of components mapping isomorphically to 
$U$, is the divisor of singularities of $\varphi _U$. The subvariety $B_U\subset X_U$ is the family of 
base loci for the Lefschetz pencils.

\begin{lemma}\label{relativeLefschetz}
There is a decomposition of local systems over $U$,
$$
R^i\alpha _*\Q = R^0q_*(R^i\varphi _{*}\Q) \oplus  R^1q_*(R^{i-1}\varphi _{*}\Q) \oplus R^2q_*(R^{i-2}\varphi _{*}\Q).
$$
\end{lemma}
\begin{proof}
See \cite[p.112]{De2}.
\end{proof}

We need a de Rham version of this statement. It would be good to have a direct analogue, however the singularities in
the Lefschetz fibration are not normal crossings singularities. Thus, we consider an open subset $W$ obtained by removing
the singular fibers. In the start of the next section we will give a series of reductions which basically say that
it suffices to look at this open set. 

It should be possible to obtain a direct de Rham version of Lemma \ref{relativeLefschetz} by applying the general
theory of $\cD$-modules, see \cite{ConsaniKim}.

Consider $\cO _{Z_U}$ the trivial vector bundle with standard connection $d$. This corresponds to the
constant local system.

Taking the relative de Rham cohomology we obtain a 
vector bundle together with its Gauss-Manin connection
$$
\cH ^i_{DR}(Z_U/U, (\cO _{Z_U}, d)) := \R ^i(f\circ \alpha )_{\ast}DR (Z_U/U, \cO _{Z_U}, d)
$$
on $U$. Let $W:= P_U - K_U$ be the complement of the locus of singularities of the family of 
Lefschetz pencils. Let $Z_W:= \varphi ^{-1}(W)\subset Z$. The map $\varphi$ is smooth over $W$.  

Similarly using the morphism $\varphi$, define $(E^i_W,\nabla _W )$ to be the vector bundle with Gauss-Manin connection
$$
E^i_W:= \cH ^i_{DR}(Z_W/W, (\cO _{Z_W}, d)):= \R ^i\varphi _{\ast}DR (Z_W/W, \cO _{Z_W}, d)
$$
As usual this extends over $P$ to a logarithmic connection with rational residues \cite{Steen}. 
Let $(E^i,\nabla )$ denote the associated parabolic
bundle on $P$, and let $(E^i_U, \nabla _U)$ denote the restriction to $P_U$ which is also the associated parabolic
bundle on $P_U$. Note that $(E^i,\nabla )$ has singularities along the divisor $D:=K\cup q^{-1}(J)$ so 
$(E^i_U, \nabla _U)$ has singularities along $K_U$.

For any multi-index $\alpha$ for the divisor $D$ we obtain a vector
bundle with logarithmic connection  $(E^i_{\alpha}, \nabla )$ on $P$ with singularities along $D$. Let
$(E^i_{\alpha ,U}, \nabla _U)$ be its restriction to $P_U$. This is the same as the corresponding bundle associated
to the parabolic bundle $E^i_U$ for multi-index $\alpha _U$ which contains only the indices for the components
of the divisor $K_U$.

In what follows we will also consider the morphism $Z_W\rightarrow U$ which doesn't have projective fibers.
We can still define the relative de Rham cohomology $\cH ^{m}_{DR}(Z_W/U, (\cO _{Z}, d))$ which is 
again a vector bundle with integrable Gauss-Manin connection on $U$. This has regular singularities and 
corresponds to the local system
$R^m(\varphi _W\circ q_W)_{\ast}\comx _{Z_W}$ on $U$. For our present purposes this can be considered as a matter of notation:
the $\cH ^{m}_{DR}(Z_W/U, (\cO _{Z}, d))$ can be defined as the unique bundle with regular singular integrable
connection on $U$ corresponding to the local system $R^m(\varphi _W\circ q_W)_{\ast}\comx _{Z_W}$, so we don't need
to explain in detail how to construct it which would involve resolving the double point singularities in the Lefschetz
fibration and then taking relative logarithmic de Rham cohomology. 

\begin{lemma}
\label{lefschetzlerayW}
Assume that  $\alpha _i >0$  for indices corresponding to components of the divisor $K$. 
The Leray spectral sequence for $Z_W\rightarrow W\rightarrow U$ is
$$
E^{i,j}_2 =\cH ^i_{DR}(P_U/U, E^j_{\alpha ,U}, \nabla ) \Rightarrow \cH ^{i+j}_{DR}(Z_W/U, (\cO _{Z}, d)). 
$$
This is a spectral sequence
in the category of vector bundles over $U$ with integrable connection having regular singularities at 
the complementary divisor $J$. 
\end{lemma}
\begin{proof}
The category of vector bundles over $U$ with integrable connection having regular singularities at 
the complementary divisor $J$ is equivalent to the category of local systems over $U$, by 
\cite{De}. Therefore it suffices to have the spectral sequence for the associated local systems.

The Leray spectral sequence for $Z_W\stackrel{\varphi _W}{\rightarrow} W\stackrel{q_W}{\rightarrow} U$ is
$$
E^{i,j}_2 =R^iq_{W,\ast}(R^j \varphi _{W,\ast} \comx _W) \Rightarrow R^{i+j}(q_W\circ \varphi _W)_{\ast}\comx _W.
$$
By Corollary \ref{logcohpar} in each fiber, we have that 
$$
\cH ^i_{DR}(P_U/U, E^k_{\alpha ,U}, \nabla ) 
$$
is the vector bundle with regular singular integrable connection corresponding to the local system 
$$R^iq_{W,\ast}(R^j \varphi _{W,\ast} \comx _W).$$ Note that the Gauss-Manin connection on 
$\cH ^i_{DR}(P_U/U, E^j_{\alpha ,U}, \nabla )$ is known to have regular singularities, in fact we shall use its
extension as a logarithmic connection given by the relative logarithmic de Rham cohomology over $S$. 

On the other hand, 
$$
\cH ^{i+j}_{DR}(Z_W/U, (\cO _{Z}, d))
$$
is the vector bundle with integrable connection corresponding to the local system 
$R^{i+j}(q_W\circ \varphi _W)_{\ast}\comx _W$. If we use the paragraph before the statement of the present lemma
then this is by definition, even. 

In view of the equivalence of categories, we get the desired spectral sequence. 
\end{proof}

Let 
$$
\overline{\cH} ^i_{DR}(P_U/U, E^j_{\alpha,U} , \nabla )
$$
and 
$$
\overline{\cH} ^{i+j}_{DR}(Z_W/U, (\cO _{Z}, d))
$$
denote the parabolic bundles on $S$ associated to these vector bundles with logarithmic connection on $U$.

\begin{corollary}
\label{lefschetzLerayFormula}
Assume that  $\alpha _i >0$  for indices corresponding to components of the divisor $K$. 
In the Chow group of $S$ tensored with $\Q$ we have
$$
\sum _m (-1)^m {\rm ch}^{\rm par}\overline{\cH} ^{m}_{DR}(Z_W/U, (\cO _{Z}, d)) \,=\,
\sum _{i,j}(-1)^{i+j}{\rm ch}^{\rm par} \overline{\cH} ^i_{DR}(P_U/U, E^j_{\alpha , U}, \nabla ).
$$
\end{corollary}
\begin{proof}
Formation of the associated parabolic bundle commutes with exact sequences by Lemma \ref{le.-parabolic2}.
The spectral sequence of Lemma \ref{lefschetzlerayW} thus gives the result. 
\end{proof}

\section{Main Theorem}

\subsection{The main statement}

We now give the statement of our main theorem in the form suitable for our inductive argument.

Suppose $S$ is a smooth projective variety and $A\subset S$ is a normal crossings divisor with smooth components.
Let $U:= S-A$. Suppose $f:X\rightarrow S$ is a projective morphism which is smooth over $U$, and let $X_U:= X\times _SU$. 
Let
$$
\cH ^i _{DR}(X_U/U):= \R ^i f_{\ast} (\Omega ^{\cdot}_{X_U/U}, d)
$$
denote the
bundle of $i$-th relative de Rham cohomology bundle on $U$. It has a Gauss-Manin connection $\nabla $
which is well-known
to have regular singularities and rational residues (\cite[Proposition 2.20]{Steen}). Let
$\overline{\cH}^i_{DR}(X_U/U)$ denote
the parabolic bundle associated to $(\cH ^i_{DR}(X_U/U), \nabla )$ on $U$.

We extend this notation to the case when $f:X\rightarrow S$ is not generically smooth. There is still an
open set $U\subset S$ over which the map is topologically a fibration for the usual topology, and we assume that
the complement has normal crossings. In this case, define $\cH ^i _{DR}(X_U/U)$ to be the unique vector bundle
with regular singular connection still denoted $\nabla $ which corresponds to the local system
$R^if_{U,\ast}\comx _{X_U}$ over $U$. Again let $\overline{\cH}^i_{DR}(X_U/U)$ denote
the associated parabolic bundle on $S$. 

We state a general form of our main theorem. It will be useful to have this generality in the inductive proof.

\begin{theorem}\label{th.-deRham}
Fix $n$. Suppose $X,A,U,Y,f$ are as above, including the case where $f$ is  projective but not generically smooth. 
Suppose that  $f: X\rightarrow S$ has relative dimension $\leq n$.
Then the alternating sum of Chern characters of the parabolic extensions to $S$
$$
\sum _{i=0}^{2n} (-1)^i{\rm ch}^{\rm par}(\overline{\cH}^i(X_U/U))
$$
lies in $CH^0(S)_\Q$ or equivalently the pieces in all of the positive-codimension Chow groups with rational coefficients
vanish.
\end{theorem}

The proof will be by induction on $n$ and using Lefschetz pencils. \footnote{Bloch and Esnault \cite{BE2} used a similar approach earlier to prove a Riemann--Roch statement for Chern--Simons class which essentially required to consider the value at the generic point of the base $S$.}  

\subsection{Preliminary reductions}
We start with some reductions. The first one 
was discussed in \S 3.5, Lemma \ref{modif-pullback} and Corollary \ref{modif}. 

\begin{lemma}
\label{reduction1}
In proving the theorem for $f:X\rightarrow S$, we can modify $S$ by a generically finite morphism $S'\rightarrow S$ and it
suffices to prove the theorem for $X'\rightarrow S'$.
\end{lemma}

\begin{lemma}
\label{reduction2}
In order to prove the theorem for $n$ and for morphisms $X\rightarrow S$ which may not be smooth, 
it suffices to prove the theorem for any $Y_U\rightarrow U$ a 
smooth projective family of relative
dimension $\leq n$.
\end{lemma}
\begin{proof}
We prove the present lemma also by induction on $n$. Therefore, we may consider a family $f:X\rightarrow S$
of relative dimension $n$
and assume that the theorem is known for arbitrary (not necessarily smooth) projective families of relative dimension
$<n$. 

After possibly making a modification of the base $S$ as in Lemma \ref{reduction1}, specially to decrease the
size of the open set $U$ and resolve singularities of the complementary divisor, we can assume that we have
an open set $W\subset X$ such that $f:W\rightarrow U$ is smooth and topologically a fibration. Let $T:= X-W$
and assume that $(X_U,T_U)\rightarrow U$ is topologically a fibration. 

Choose another compactification $W\subset Y$ with a morphism $g:Y\rightarrow U$ such that $g$ is smooth and projective.
Let $V:= Y-W$ be the complementary divisor which we assume to have relative normal crossings over $U$. 

Fix a point $u\in U$. We have pairs $(X_u, T_u)$ and $(Y_u, V_u)$. In both of these, the smooth complementary open
set is $W_u$. The cohomology of either of these pairs is equal to the compactly supported cohomology of $W_u$. 
Poincar\'e duality says 
$$
H^i((X_u,T_u),\comx ) \cong H^{2n-i}(W_u, \comx ) \cong H^i((Y_u,V_u),\comx ).
$$
On the other hand, we have a long exact sequence relating $H^{\cdot}(X_u, \comx )$,
$H^{\cdot}(T_u, \comx )$, and $H^i((X_u,T_u),\comx )$. Similarly we have a long exact sequence
relating $H^{\cdot}(Y_u, \comx )$,
$H^{\cdot}(V_u, \comx )$, and $H^i((Y_u,V_u),\comx )$.  As $u$ varies in $U$, the Poincar\'e duality isomorphisms
are isomorphisms of local systems over $U$, and the long exact sequences are long exact sequences in the category
of local systems over $U$. In view of the fact that the associated parabolic bundle commutes with isomorphisms and
with taking exact
sequences, we obtain a long exact sequence relating
$$
\overline{\cH}^{\cdot} _{DR}(X_U/U), \;\;\; \overline{\cH}^{\cdot} _{DR}(T_U/U), \;\; \mbox{and} \;
\overline{\cH}^{2n-\cdot} _{DR}(W/U),
$$
and a long exact sequence relating
$$
\overline{\cH}^{\cdot} _{DR}(Y_U/U), \;\;\; \overline{\cH}^{\cdot} _{DR}(V_U/U), \;\; \mbox{and} \;
\overline{\cH}^{2n-\cdot} _{DR}(W/U).
$$
The shift of indices induced by Poincar\'e duality is $2n$, an even number, so it doesn't affect the sign
of the alternating sum. Thus we obtain
$$
\sum _{i=0}^{2n} (-1)^i{\rm ch}(\overline{\cH}^i(X_U/U)) - 
\sum _{i=0}^{2n-2} (-1)^i{\rm ch}(\overline{\cH}^i(T_U/U)) = 
\sum _{i=0}^{2n} (-1)^i{\rm ch}(\overline{\cH}^i(W/U))
$$
and
$$
\sum _{i=0}^{2n} (-1)^i{\rm ch}(\overline{\cH}^i(Y_U/U)) - 
\sum _{i=0}^{2n-2} (-1)^i{\rm ch}(\overline{\cH}^i(V_U/U)) = 
\sum _{i=0}^{2n} (-1)^i{\rm ch}(\overline{\cH}^i(W/U)).
$$
In particular, 
by the induction hypothesis we know the statement of the theorem for $T\rightarrow S$ and $V\rightarrow S$
because these have relative dimension $\leq n-1$. On the other hand we also know the statement of the theorem
for $Y\rightarrow S$, because this is the hypothesis of our present reduction lemma. Putting these together we obtain
the statement of the theorem for $X\rightarrow S$.
\end{proof}

\begin{lemma}
The statement of the theorem is true in relative dimension $n=0$. 
\end{lemma}
\begin{proof}
Indeed, in this case by making a generically finite modification of the base, we can reduce to the case when
the morphism $X_U\rightarrow U$ is a disjoint union of sections isomorphic to $U$. Thus the relative cohomology
sheaf, nontrivial only in degree zero, is a trivial bundle.
\end{proof}

{\bf Remark:} It would be good to have the statement of Theorem \ref{th.-deRham} for families which are not necessarily projective.
Our method of proof gives the statement for smooth quasiprojective families. The reduction from the general case to 
these cases would seem to require a more detailed utilisation of Grothendieck-Verdier 
duality theory than we are prepared to do here. 

\subsection{Lefschetz pencil and reduction to the case of the open set $Z_W/U$}

We now assume given a projective family $f:X\rightarrow S$ of relative dimension $n$ 
such that $X_U$ is smooth over $U$,
and assume the theorem is known for relative dimension $\leq n-1$. After modification of the base as per 
Lemma \ref{modif-pullback},
arrange to have a family of Lefschetz pencils as in \S 4.2. Keep the same notations $Z,P,W$ from there.

We assume that we have modified $S$ enough so that the monodromies of the local systems which we encounter
on $U$ or $W$ are unipotent
around components of $J$ or vertical components of $f^{-1}(J)$. In this case, the extended bundles from $U$ to
$S$ are actual bundles rather than parabolic bundles. We may use the same notation ${\rm ch}$ for the parabolic Chern 
character of parabolic bundles as that
used for regular bundles, but where possible we indicate when a nontrivial parabolic structure might be involved by the notation
${\rm ch}^{\rm par}$. 

By our inductive hypothesis, we know the statement of the theorem for the morphisms 
$$
\varphi : Z\rightarrow P, \;\;\; q:P\rightarrow S.
$$
Let $(E^i,\nabla )$ denote the parabolic bundles considered in \S 4.2. By the inductive hypothesis for the map $\varphi$
we have
$$
\sum _{i=0}^{2n} (-1)^i{\rm ch}(E^i)  \in CH^0(P)_{\Q} .
$$
We are assuming from the previous paragraph that the monodromies of the local system $E^{\nabla}$ are
unipotent around components of $f^{-1}(J)$, so the parabolic structure on $E$ along these components of the divisor
$D$ is a trivial parabolic structure on a usual bundle.  It is only in the case of a Lefschetz pencil of even
fiber dimension (that is, when $n-1$ is even) that $E$ has nontrivial parabolic structure along the horizontal
components $K$. 

We show in the next two lemmas that the main step in the proof will be to treat the case $Z_W/U$.

\begin{lemma}
\label{reductionZ_UtoZ_W}
Suppose we know the statement of the theorem for $Z_W/U$, in other words suppose we know that
$$
\sum _{i=0}^{2n} (-1)^i{\rm ch}(\overline{\cH}^i(Z_W/U)) \in CH^0(S)_{\Q} .
$$
Then we can conclude the statement of the theorem for $Z_U/U$.
\end{lemma}
\begin{proof}
Let $Z_{K,U}:= Z_U-Z_W$. It is equal to the preimage $\varphi _U ^{-1}(K_U)$, so it is
a closed subset of $Z_U$. 

Proceed as in the proof of Lemma \ref{reduction2}. The pair $(Z_U, Z_{K,U})$ is a fibration over $U$.
For each $u\in U$ the relative cohomology $H^{\cdot}((Z_u, Z_{K,u}),\comx )$ fits into a long exact sequence
with $H^{\cdot}(Z_u, \comx )$ and $H^{\cdot}(Z_{K,u}, \comx )$. On the other hand, 
$H^{\cdot}((Z_u, Z_{K,u}),\comx )$ is isomorphic to the compactly supported cohomology of $W_u$ 
which is Poincar\'e dual to the cohomology of $W_u$ with $i\mapsto 2n-i$. This shift doesn't affect the signs in
our alternating sums. 
As before, we obtain long exact sequences of local systems on $U$ and then of associated parabolic bundles on $S$.
Thus, the hypothesis of the present lemma implies that 
$$
\sum _{i=0}^{2n} (-1)^i{\rm ch}(\overline{\cH}^i(Z_{U}/U)) -
\sum _{i=0}^{2n} (-1)^i{\rm ch}(\overline{\cH}^i(Z_{K,U}/U)))\in CH^0(S)_{\Q}.
$$
On the other hand the morphism $Z_{K,U}\rightarrow U$ is projective of relative dimension $n-1$. 
It is not smooth, but we have taken care to include the non-smooth case in our inductive statement. Thus,
by the inductive assumption we know that 
$$
\sum _{i=0}^{2n} (-1)^i{\rm ch}(\overline{\cH}^i(Z_{K,U}/U)))\in CH^0(S)_{\Q} .
$$
We obtain the same conclusion for $Z_U/U$.
\end{proof}

\begin{lemma}
\label{reductionX_UtoZ_U}
Suppose we know the statement of the theorem for $Z_U/U$, in other words suppose we know that
$$
\sum _{i=0}^{2n} (-1)^i{\rm ch}(\overline{\cH}^i(Z_U/U)) \in CH^0(S)_{\Q} .
$$
Then we can conclude the statement of the theorem for $X_U/U$.
\end{lemma}
\begin{proof}
Recall that $B\subset X$ is the family of base loci for the Lefschetz pencil, and $Z$ is obtained from $X$ by 
blowing up along $B$, at least in the part lying over $U$. Let $\tilde{B}$ denote
the inverse image of $B$ in $Z$. The relative dimension of $B_U/U$ is $n-2$ and the relative dimension of 
$\tilde{B}_U/U$ is $n-1$. Thus we can assume that we know the statement of the theorem for these families. 
Furthermore $Z_U -\tilde{B}_U \cong X_U - B_U$ and this is a smooth quasiprojective variety with topological
fibration to $U$. Use the same argument as in Lemma \ref{reduction2} applying the long exact sequences for
pairs and Poincar\'e duality on $(X_U -B_U)/U$. By this argument applied once to each side, 
if we know the statement of the theorem for
$Z_U/U$ then we get it for $(Z_U -\tilde{B}_U)/U = (X_U -B_U)/U$ and hence for $X_U/U$.
\end{proof}

\begin{corollary}
\label{reductionX_UtoZ_W}
Suppose we know the statement of the theorem for $Z_W/U$, 
then we can conclude the statement of the theorem for $X_U/U$.
\end{corollary}
\begin{proof}
Put together the previous two lemmas.
\end{proof}

\subsection{The proof for $Z_W/U$}

Up to now we have reduced to the problem of proving the statement of the theorem for 
the family of open smooth varieties $Z_W/U$. Here is where we use the Lefschetz pencil
$Z_W\stackrel{\varphi}{\rightarrow} W \stackrel{q}{\rightarrow} U$ as in \S 4. 

Recall that we defined
$$
E^i_W:= \cH ^i_{DR}(Z_W/W, (\cO _{Z_W}, d)):= \R ^i\varphi _{\ast}DR(Z_W/W, \cO _{Z_W}, d)
$$
with its Gauss-Manin connection $\nabla _W$ on $W$. Then  $(E^i,\nabla )$ was the associated parabolic
bundle on $P$, with $(E^i_U, \nabla _U)$  the restriction to $P_U$.
For any multi-index $\alpha$ for the divisor $D$ we have $(E^i_{\alpha}, \nabla )$ 
which is a vector bundle
on $P$ with logarithmic connection having singularities along $D$, and
$(E^i_{\alpha ,U}, \nabla _U)$ is the restriction to $P_U$. 

We assume that $\alpha _i >0$ for components $D_i$ appearing in $K$, that is components which surject to $S$,
whereas we assume from now on that $\alpha _i = 0$ for vertical components $D_i$, that is the components
which map to components of $J$. 

By Corollary  \ref{lefschetzLerayFormula} we have
$$
\sum _m (-1)^m {\rm ch}\overline{\cH} ^{m}_{DR}(Z_W/U) \, = \, 
\sum _{i,j}(-1)^{i+j}{\rm ch} \overline{\cH} ^i_{DR}(P_U/U, E^j_{\alpha , U}, \nabla ).
$$
On the left in this formula is the quantity which we are trying to show lies in  $CH^0(S)_{\Q}$.

Recall that the quantity on the right is defined to be the parabolic bundle on $S$ associated to
the vector bundle $\cH ^i_{DR}(P_U/U, E^j_{\alpha , U}, \nabla )$
with regular singular Gauss-Manin connection over $U$.

\begin{lemma}
\label{applySteenbrink}
The extended parabolic bundle $\overline{\cH} ^i_{DR}(P_U/U, E^j_{\alpha , U}, \nabla )$
is a vector bundle with trivial parabolic structure, and it corresponds to the vector bundle
$\cH ^i_{DR}(P/S, E^j_{\alpha }, \nabla )$ obtained by taking the relative logarithmic de Rham cohomology over $S$. 
\end{lemma}
\begin{proof}
This is the contents of our analogue of Steenbrink's theorem, Theorem \ref{th.-steen} (see \S6). Note that the fact that
we have chosen $\alpha _i = 0$ for vertical components of the divisor $D_i$ means that the logarithmic 
connection $\nabla $ on $E^i_{\alpha}$ has nilpotent residues along the vertical components $D_i$. Recall for 
this that we are assuming that we have made sufficient modification $S'\rightarrow S$ so that the monodromy
transformations of $(E,\nabla )$ around vertical components of the divisor $D$ are unipotent. 
\end{proof}

\begin{corollary}
\label{lastreduction}
For the proof of Theorem \ref{th.-deRham} it suffices to show that
$$
\sum _{i,j}(-1)^{i+j}{\rm ch}\cH ^i_{DR}(P/S, E^j_{\alpha }, \nabla ) \in CH^0(S)_{\Q} .
$$
\end{corollary}
\eop

We now proceed with the proof of the formula to which we have reduced in the previous corollary.

Use the inductive statement for the family $\varphi : Z\rightarrow P$. Note that this has relative dimension $n-1$. 
The parabolic bundle on $P$ associated to the higher direct image local system on the open set $W$
is $E^i$.
Therefore our inductive hypothesis for $Z/P$ says that
$$
\sum _i (-1)^i {\rm ch}^{\rm par}(E^i) \in CH^0(P)_{\Q} .
$$

Next observe that by Lemma \ref{diffpar}, the difference between the parabolic Chern character of the parabolic bundle,
and the Chern character of any component bundle, comes from the divisor $D$. In fact, due to our assumption about
the monodromy, $E^i$ is just a usual bundle along the vertical components of the divisor. Thus, for
multi-indices $\alpha$ such that $\alpha _i = 0$ along vertical components, the difference comes
from the divisor $K$:
$$
{\rm ch}^{\rm par}(E^i)-{\rm ch}(E^i_{\alpha} ) \in {\rm Image}(CH^{\cdot}(K)_{\Q} \lrar CH^{\cdot}(P)_{\Q}) .
$$
To be more precise this means that the difference is in the subspace of $CH^{\cdot}(P)_{\Q} $ spanned by
Chern characters of sheaves concentrated on $K$. Thus we get
$$
\sum _i (-1)^i{\rm ch}(E^i_{\alpha} )\in CH^{\cdot}(K)_{\Q}  + CH^0(P)_{\Q} .
$$
Write
$$
\sum _i (-1)^i{\rm ch}(E^i_{\alpha} ) = {\bf k} + {\bf r}
$$
with ${\bf k}\in CH^{\cdot}(K)_{\Q}$ and ${\bf r}\in CH^0(P)_{\Q}$. Note that ${\bf r}$ is just an 
integer, the alternating sum of the ranks of the bundles.

Now look at the relative logarithmic de Rham complex for $E^i_{\alpha}$:
$$
E^i_{\alpha}\rightarrow E^i_{\alpha}\otimes _{\cO _P}\Omega ^1_{P/S}(\log D).
$$
The Hodge to de Rham spectral sequence and multiplicativity for tensor products and additivity for exact sequences
of the Chern character, gives the formula
$$
\sum _{i,j}(-1)^{i+j}{\rm ch}\cH ^i_{DR}(P/S, E^j_{\alpha }, \nabla ) \, = \,
$$
$$
\chi _{P/S}\left( \left[ \sum _i (-1)^i{\rm ch}(E^i_{\alpha} )\right] 
\cdot \left[\sum _{j=0,1} (-1)^j {\rm ch} \, \Omega ^j_{P/S}(\log D)\right] \right) .
$$
Recall here that $\chi _{P/S}$ represents the function (see Proposition \ref{pr.-GRR}), entering into the GRR formula for $P/S$ calculating
the relative Euler characteristic in $CH^{\cdot}(S)_{\Q}$ for an element of $CH^{\cdot}(P)_{\Q}$.

This formula decomposes into two pieces according to the decomposition ${\bf k} + {\bf r}$ above:
\begin{eqnarray*}
\sum _{i,j}(-1)^{i+j}{\rm ch}\cH ^i_{DR}(P/S, E^j_{\alpha }, \nabla ) & = &
\chi  _{P/S}\left( {\bf k}
\cdot \left[\sum _{j=0,1} (-1)^j {\rm ch} \, \Omega ^j_{P/S}(\log D)\right] \right) \\
&& +\,\, 
\chi  _{P/S}\left( {\bf r}
\cdot \left[\sum _{j=0,1} (-1)^j {\rm ch} \, \Omega ^j_{P/S}(\log D)\right] \right) .
\end{eqnarray*}

\begin{lemma}
\label{firstpiece}
In the above decomposition, the first piece vanishes: 
$$
\chi  _{P/S}\left( {\bf k}
\cdot \left[\sum _{j=0,1} (-1)^j {\rm ch} \, \Omega ^j_{P/S}(\log D)\right] \right) =0
$$
\end{lemma}
\begin{proof}
We claim that the argument inside the function $\chi  _{P/S}$ is already zero. 
Indeed, ${\bf k}$ is a sum of Chern characters of sheaves supported on $K$,
but on the other hand the divisor $K$ is a union of components which appear in $D$ and only intersect
single vertical components of the divisor $D$, because $(P,K)\rightarrow S$ is a semistable family of pointed
curves. Thus we have a residue isomorphism
$$
\Omega ^1_{P/S}(\log D) |_ K \cong \cO _K = \Omega ^0_{P/S}(\log D) |_K,
$$
so the difference $\sum _{j=0,1} (-1)^j {\rm ch} \, \Omega ^j_{P/S}(\log D)$ restricts to zero on $K$. When
we multiply it by ${\bf k}$ we get zero. 
\end{proof}

After this we are left only with the second piece:
$$
\chi  _{P/S}\left( {\bf r}
\cdot \left[\sum _{j=0,1} (-1)^j {\rm ch} \, \Omega ^j_{P/S}(\log D)\right] \right) .
$$
The class ${\bf r}$ is just an integer considered in $CH^0(P)_{\Q}$. Thus, the following lemma will complete
the proof.

\begin{lemma}
In our situation, 
$$
\chi  _{P/S}\left( \sum _{j=0,1} (-1)^j {\rm ch} \, \Omega ^j_{P/S}(\log D) \right) \in CH^0(S)_{\Q} .
$$
\end{lemma}
\begin{proof}
This quantity is equal to 
$$
\sum _{i}(-1)^{i+j}{\rm ch}\cH ^i_{DR}(P/S, (\cO _P, d) ).
$$
By the same argument as above using Steenbrink's theorem for the fibration $q:P\rightarrow S$, the
terms $\cH ^i_{DR}(P/S, (\cO _P, d) )$ are the vector bundles over $S$ assocated by the Deligne canonical extension
to the local systems 
$R^iq_{U,\ast}\comx _{P_U}$ defined on $U$. Since $P/S$ is semistable the monodromy here is unipotent. 
However in this case $P_U\rightarrow U$ is a family of smooth rational curves. Thus the direct image local systems
are $\comx _U$ for $i=0,2$ and $0$ for $i=1$. Their canonical extensions are trivial bundles on $S$ so
$$
\sum _{i}(-1)^{i+j}{\rm ch}\cH ^i_{DR}(P/S, (\cO _P, d) ) \in CH^0(S)_{\Q} .
$$
\end{proof}

Multiplied by the integer ${\bf r}$, the result of this lemma shows that the second term in the
previous decomposition lies in $CH^0(S)_{\Q} $. As we saw in Lemma \ref{firstpiece} that the other piece
vanishes, we obtain:

\begin{corollary}
\label{endproof}
$$
\sum _{i,j}(-1)^{i+j}{\rm ch}\cH ^i_{DR}(P/S, E^j_{\alpha }, \nabla )\in CH^0(S)_{\Q} .
$$
\end{corollary}
\eop

In view of Corollary \ref{lastreduction}, this completes the proof of Theorem \ref{th.-deRham}.
\eop

\subsection{Family of projective surfaces}
 
Consider a family of projective surfaces $$\pi:X\lrar S$$ where $S$ is smooth. 
Then, with notations as in Theorem \ref{th.-deRham}, we have

\begin{proposition}\label{pr.-surfaces}
The parabolic Chern character of the parabolic bundles associated to weight $i$ de Rham bundles satisfies
$$ch(\ov\cH^i(X/S))\,\in\, CH^0(S)_\Q$$
for each $i\geq 0$.
\end{proposition}

We first consider the following case:

\begin{lemma}\label{le.-smooth}
Suppose $\pi:X\lrar S$ is a generically smooth morphism of relative dimension $2$. Then
$$ch(\ov\cH^i(X/S))\,\in\, CH^0(S)_\Q$$
for each $i\geq 0$.
\end{lemma}
\begin{proof}
Firstly, the Hard Lefschetz theorem and \cite[Theorem 1.1]{Es-Vi2} imply that
$$ch(\ov\cH^1(X/S)),\, ch(\ov\cH^3(X/S)) \,\in\, CH^0(S)_\Q.$$
Hence by Theorem \ref{th.-deRham}, we conclude that
$$ch(\ov\cH^2(X/S))\,\in\, CH^0(S)_\Q.$$
\end{proof}

Now suppose $\pi:X\lrar S$ be any projective morphism.
Let $\m{Spec}\,K$ be the generic point of $S$ and $X_K:=X\times _\comx \m{Spec}\,K$. Since the field $K$ is of 
characteristic zero, by Hironaka's theorem \cite{Hi}, there is a resolution of singularities $$f:X'_K\lrar X_K$$ 
defined over the field $K$ such that $f$ is a sequence of blow-ups and $X'_K$ is nonsingular.

This implies that there is a non-empty open subset $U\subset S$ and a commutative diagram:
\begin{eqnarray*}
X'_U &\sta{f}{\lrar} & X_U \\
\downarrow \pi'&  & \downarrow\pi\\ 
U \,\,\,    &=& U
\end{eqnarray*}
such that $\pi'$ is a smooth morphism and the variety $X'_U\lrar U$ is a 
fibrewise resolution of singularities of $X_U\lrar U$. By Lemma \ref{modif-pullback}, 
we assume that $D:=S-U$ is a normal crossing divisor. 

Since a fibre $X_s$ of $\pi$ is a singular surface, the cohomology of $X_s$ carries a mixed Hodge structure. 
In particular, over the open subset of $S$ where $\pi$ is a topological fibration 
(which again by Lemma \ref{modif-pullback} can be assumed to be $U$), we obtain a filtration 
of local systems with associated-graded
\begin{equation}\label{eq.-local}
Gr _{\cdot}R^i\pi_*\Q\,=\,\oplus_{k=0}^i Gr_kR^i\pi_*\Q.
\end{equation}

Here $Gr_kR^i\pi_*\Q$ is the local system corresponding to the weight $k$ graded 
piece of the weight filtration $\{W_k\}$on the cohomology $H^i(X_s,\Q)$. 
More precisely, $H^i(X_s,\Q)$ carries a weight filtration 
$$(0)=W_{-1}\subset W_0\subset...\subset W_i=H^i(X_s,\Q)$$
such that $Gr_kH^i(X_s,\Q)$ carries a polarised pure Hodge structure of weight $k$.
In other words, $Gr_kR^i\pi_*\Q$ is the local system associated to a family of polarised 
pure Hodge structures of weight $k$. We refer to \cite{De3} for the details.

Also, notice that the morphism $f^*_s:W_i= H^i(X_s,\Q)\lrar H^i(X'_s,\Q)$ is a morphism of 
mixed Hodge structures whose kernel is $W_{i-1}$. In particular, the $i$-th graded piece 
$Gr_iH^i(X_s,\Q)$ is a polarised pure sub-Hodge structure of weight $i$ of $H^i(X'_s,\Q)$. 
When $i=2$, the complementary sub-Hodge structure is generated by the algebraic classes since
$X'_s\lrar X_s$ is a sequence of blow-ups and normalizations
for $s\in U$.  

These statements put together over $U$ corresponds to saying that 

\begin{lemma}\label{le.-triv}
There is a decomposition of local systems 
$$
R^2\pi'_*\Q\,=\,Gr_2R^2\pi_*\Q \oplus T
$$
where $T$ is a trivial local system.
\end{lemma}
\begin{proof}
The kernel of the map on $H_2$ of a surface corresponding either to blowing up a point, or to normalization, is generated
by algebraic cycles. Therefore the kernel of $H_2(X'_s)\rightarrow H_ 2(X_s)$ is generated by algebraic cycles.
This gives an exact sequence
$$
0 \lrar \,Gr_2R^2\pi_*\Q \lrar R^2\pi'_*\Q \lrar T \lrar 0
$$
where $T$ corresponds to a local system of finite monodromy. After going to a finite cover of $U$ we may
assume $T$ is trivial. The cup product on $R^2\pi'_*\Q$ is a nondegenerate form which is nondegenerate on
$Gr_2R^2\pi_*\Q$ so we get an orthogonal splitting of the exact sequence. 
\end{proof}

Consider the parabolic bundles $\ov\cG_k(\cH^i)$ associated to the graded pieces $Gr_kR^i\pi_*\Q$, for each $i,k$, on $S$.

Then we have
\begin{lemma}\label{le.-graded}
The parabolic Chern character of the above parabolic bundles satisfies
$$ch(\ov\cG_k(\cH^i))\,\in\, CH^0(S)_\Q$$
for $i,k\leq 2$.
\end{lemma}
\begin{proof}
By Lemma \ref{le.-smooth}, we have $ch(\ov\cH^2(X'/S))\in CH^0(S)_\Q$. Hence by Lemma \ref{le.-triv}, we conclude that
$$ch(\ov\cG_2(\cH^2))\,\in\, CH^0(S)_\Q.$$
Since the parabolic bundle $\ov\cG_0(\cH^i)$ for any $i$, is a trivial bundle and $\ov\cG_1(\cH^i)$ corresponds
to a family of polarised pure Hodge structures of weight one, by \cite[Theorem 1.1]{Es-Vi2}, we conclude  
that
$$ch(\ov\cG_1(\cH^2))\,\in\, CH^0(S)_\Q.$$
\end{proof}

Using \eqref{eq.-local} and Lemma \ref{le.-parabolic2}, we conclude that
$\ov\cH^2(X/S)$ (resp. $\ov\cH^1(X/S)$) has a filtration whose associated-graded is $\oplus_{k=0}^2 \ov\cG_k(\cH^2)$
(resp. $\oplus_{k=0}^1 \ov\cG_k(\cH^1)$).
 
By Lemma \ref{le.-graded}, we obtain 
$$ch(\ov\cH^i(X/S))\in CH^0(S)_\Q$$
for $i=0,1,2$. Now by Theorem \ref{th.-deRham}, we obtain
$$ch(\ov\cH^3(X/S))\in CH^0(S)_\Q.$$

This completes the proof of Proposition \ref{pr.-surfaces}.

\section{Appendix: an analogue of Steenbrink's theorem}

Suppose $f:P\rightarrow S$ is a semistable morphism of relative dimension one, with respect to a divisor $J$ on $S$ and 
its pullback
$f^{-1}J$ on $P$, both of which have normal crossings. Assume for simplicity that the total space of $P$ is smooth
(which is possible by \cite[Proposition 3.6]{deJong}).

Semistability implies that the components of
$f^{-1}J$ occur with multiplicity one, and that $f$ is flat. Let $U:=S-J$ and $V:= P - f^{-1}J = f^{-1}U$ be
open sets in $S$ and $P$ respectively. Suppose $K\subset P$ is another divisor, isomorphic to a disjoint union of copies
of $S$ on which $f$ is the identity. Suppose that the components of $K$ meet $f^{-1}J$ transversally at smooth points of the
latter, so that $(P,K)\rightarrow S$ is a semistable family of pointed curves. Let $D:= K\cup f^{-1}J$.

Suppose that $(E,\nabla )$ is a logarithmic connection on $P$ with singularities along $D$ with rational residues. 
Suppose that the residual transformations of
$\nabla$ along components of $f^{-1}J$ are nilpotent. Assume also that for every component $K_i$ of the horizontal
divisor $K$ we have a weight $\alpha_i$ such that the eigenvalues of the residues of $\nabla$ along $K_i$ are contained
in $[-\alpha _i, 1-\alpha _i )$.

Let $E\otimes \Omega ^{\cdot}_{P/S}(\log D)$ be the
relative logarithmic de Rham complex (including $K$ in the logarithmic part). Set
$$
F:= \cH ^{\cdot} _{DR}(P/S, E, \nabla ) = \R f_{\ast} (E\otimes \Omega ^{\cdot}_{P/S}(\log D)).
$$
This is a complex of $\cO _S$-modules, quasi-isomorphic to a complex of vector bundles, in other words it is
a {\em perfect complex} on $S$. We can now state Steenbrink's theorem in this setting:

\begin{theorem}\label{th.-steen}
Keep the above hypothesis of a semistable family and logarithmic connection $(E,\nabla )$ such that the residues
of $\nabla$ along components of $f^{-1}J$ are nilpotent and the residues along $K_i$ have
eigenvalues in $[-\alpha _i, 1-\alpha _i )$. Then the higher direct image complex $F$ splits locally over
$S$ into a direct sum of its cohomology sheaves $H^i(F)$. Furthermore, these cohomology sheaves are locally free on
$S$. They have Gauss-Manin connections which are logarithmic along $J$, and the residues of their Gauss-Manin
connections along $J$ are nilpotent.
In particular, the cohomology sheaves $H^i(F)$ are the canonical extensions of the Gauss-Manin connections over
the open set $U$ (denoted $\overline{\cH}^i$) and
$$
{\rm ch}(F) = \sum _ i (-1)^i{\rm ch}(\overline{\cH}^i).
$$
\end{theorem}

Steenbrink proved this in the context of a geometric family \cite[Theorem 2.18]{Steen}. A weaker version of this theorem was proved earlier by N. Katz \cite{KatzICM}.  
Since then there have been a
large number of generalizations of his theorem, see
\cite{Faltings} \cite{Faltings2} \cite{Saito}  \cite{Illusie} \cite{KatoNakayama} \cite{KatoMatsubaraNakayama}
\cite{IllusieKatoNakayama} \cite{FKato} 
\cite{Tsuji} \cite{Tsuji2} \cite{Elzein} \cite{Ogus} \cite{Cailotto} among
others. The argument we sketch below is probably contained in some of these references in some way.\footnote{For example
the argument given by Illusie in \cite{Illusie} uses the condition that the higher direct image is a perfect
complex. 
Illusie doesn't treat the case of local coefficient
systems, necessary for our inductive argument. He refers to 
Faltings \cite{Faltings}, \cite{Faltings2} for treating some more general cases. The techniques we use here certainly come from this 
theory (the second author went to Faltings' course on $p$-adic Hodge theory).
Also the note of Cailotto \cite{Cailotto} seems to use an argument similar to the one we give here. 
Illusie, Kato and Nakayama touch on this result in \S 6 and Theorem 7.1 of \cite{IllusieKatoNakayama},
and refer for the case of a system of local coefficients to Kato-Matsubara-Nakayama
\cite{KatoMatsubaraNakayama}. There, this result is treated for local coefficients in a variation of Hodge structure,
which would be sufficient for our purposes since $(E,\nabla )$ which intervenes in our argument is a 
variation of Hodge structure coming from the Lefschetz pencil.}

We will indicate a proof which we hope will be enlightening. 
The first thing to note is that it suffices to prove that the dimensions of the cohomology
groups of the fibers $F_x := F \otimes _{\cO _S}{\bf k}(x)$ remain constant. Therefore it suffices to look at a general
curve going into each point, and since the semistability hypothesis is preserved under base change to a general curve
in the base (\cite{Ab-dJ}) we may assume that $S$ is a curve. Note also that the Gauss-Manin connection is globally defined
before we restrict to a curve, but it suffices to look over a curve in order to prove that the residues are nilpotent,
so for this part also it suffices to assume that $S$ is a curve.

The idea is that one can define the Gauss-Manin connection on the level of the complex $F$, Lemma \ref{actionL}
below \cite{KatzICM}.  We then observe in Proposition 
\ref{criterion} that these standard facts imply local freeness of the cohomology sheaves by a direct calculation:
any cohomology sheaves not locally free would lead to eigenvalues of the residues differing by nonzero integers which
contradicts nilpotency. This is basically the same as the main lemma in \cite{Cailotto}. 

Let $\cL$ denote the
sheaf of differential operators of order $\leq 1$, whose symbols in degree $1$ are vector fields tangent to $J$.
Since we are assuming that $S$ is a curve, an action of
$\cL$ on $F$ will contain all of the information that we need.
In order to fix the ideas, note several facts about $\cL$.   It has left and right structures of
$\cO _S$-module which are different. We have an exact sequence, compatible with both the left and right module structures:

\begin{equation}\label{eq.-c}
0\rightarrow \cO _S \rightarrow \cL \rightarrow (\Omega ^1_S(\log J))^{\ast} \rightarrow 0.
\end{equation}

In view of the assumption that $S$ is a curve, so $J$ is a collection of points, we have that
$$
(\Omega ^1_S(\log J))^{\ast} = T_S(-J)
$$
is just the sheaf of tangent vector fields vanishing at the points of $J$. Thus we can write the exact sequence as
$$
0\rightarrow \cO _S \rightarrow \cL  \stackrel{\sigma}{\rightarrow} T_S(-J) \rightarrow 0.
$$
The fiber $T_S(-D)_x$ over a point $x\in J$ is a complex line generated by the canonical element denoted
$(t\frac{\partial}{\partial t} )_x$ which is independent of the choice of coordinate $t$ at $x$. 

In the above exact sequence, 
the map $\sigma$ is the {\em symbol} of a differential operator. Using it we can describe the difference between
the left and right structures of $\cO _S$-module. For $v\in \cL $ and $a\in \cO _S$ we have
$$
v \cdot a - a \cdot v = \sigma (v)(a),
$$
this answer $\sigma (v)(a)$ means the derivative of $a$  along the vector field $\sigma (v)$, and $\sigma (v)(a)$ 
is considered as
an element of $\cO _S \subset \cL$.

Suppose $x\in J$. The fiber of $\cL$ over $x$ is the same when it is calculated on the left or the right.
Let ${\bf m}_x$ be the maximal ideal at $x$ and $\comx _x := \cO _S /{\bf m}_x$. We have
$$
\comx _x \otimes _{\cO _S} \cL = \cL \otimes _{\cO _S} \comx _x  = 
\comx \oplus \langle t\frac{\partial}{\partial t} \rangle _x .
$$
Denote this by $\cL _x$. 
This is special for points $x\in J$ because we have taken differential operators generated by vector fields which are
tangent to $J$. In particular it means that an action of $\cL$ induces an action of $\cL _x$ on the fiber over $x$. 

For general points in $S$ the quotient of $\cL$ by the maximal ideal on the left and on the right are
not canonically isomorphic, and we could no longer restrict to an action on the fiber in this way.

\begin{lemma}
\label{actionL}
Let $F:= \cH ^{\cdot} _{DR}(P/S, E, \nabla )$ be the relative logarithmic de Rham cohomology, under the hypothesis 
of Theorem \ref{th.-steen}, specifically let it be the \v{C}ech complex obtained from an affine open covering. 
There is an action of $\cL$ on $F$ given by 
$$
\cL \otimes _{\cO _S}F \rightarrow F
$$
whose restriction to the scalars $\cO _S \subset \cL$ is the identity. Furthermore if $x\in J$ then the action of
the vector $(t\frac{\partial}{\partial t})_x\in \cL _x$ on the fiber $F_x$ is nilpotent. 
\end{lemma}
\begin{proof}
Katz constructs an action of the vector field $t\frac{\partial}{\partial t}$ by derivations on the full higher direct image
complex in \cite{KatzICM}, see also \cite{KatzOda}. This is the required action of $\cL$, inducing the identity
on the scalars. Katz shows that the indicial polynomial of the action of $(t\frac{\partial}{\partial t})_x$ divides
a product of indicial polynomials of the residues of $\nabla$ on vertical components. In our case where the residues of
$\nabla$ are supposed to be nilpotent in the hypothesis of \ref{th.-steen}, we get that 
$(t\frac{\partial}{\partial t})_x\in \cL _x$ is nilpotent. \footnote{We thank H. Esnault for pointing out the reference \cite{KatzICM}
which replaces the lengthy discussion of this proof, using a truncated version of
the theory of formal categories \cite{Berthelot}, in our original version. 
One could also now invoke either the theory of $\cD$-modules,
or the {\em log-crystalline site} to get the same statement.} 
\end{proof}

Now note that $\cL$ is locally free as a
right $\cO _S$-module, so the tensor product is also the derived tensor product
$$
\cL \otimes _{\cO _S} F = \cL \otimes ^{\bbL}_{\cO _S} F.
$$
In this spirit, we have a homotopy invariance of the existence of this action.

\begin{lemma}
\label{homotopyinvariance}
Locally over $S$, suppose $F'$ is a perfect complex, that is a complex of vector bundles of finite length, 
quasiisomorphic to $F$.Then there is still an action 
$$
\cL \otimes _{\cO _S} F'\rightarrow F'
$$
with the property that for any $x\in J$, the action of $\cL _x$ on $H^i(F'_x)$ induces the identity action of the scalars
$\comx \subset \cL _x$ and a nilpotent action of the vector $(t\frac{\partial}{\partial t})_x$.
\end{lemma}
\begin{proof}
Being quasiisomorphic means that there is a chain of quasiisomorphisms going in different directions relating $F'$ 
and $F$. However, since $F'$ consists of projective objects locally in the Zariski topology of $S$, we can represent 
this by an actual morphism of complexes $F'\rightarrow F$. Here we allow ourselves to replace $S$ by a smaller neighborhood
of the point $x\in J$ we are interested in. Similarly, $\cL \otimes _{\cO _S} F'$ is a complex of vector bundles
so the map 
$$
Hom (\cL \otimes _{\cO _S} F' , F')\rightarrow Hom (\cL \otimes _{\cO _S} F' , Q^{\cdot}\otimes \langle d\log t \rangle ^{\ast})
$$
is a quasiisomorphism of complexes over $S$, and again possibly after going to an open set, our map
$$
\cL \otimes _{\cO _S} F' \rightarrow Q^{\cdot}\otimes \langle d\log t \rangle ^{\ast}
$$
lifts to a map
$$
\cL \otimes _{\cO _S} F'\rightarrow F'
$$
up to addition of a map of the form $d(\kappa )$, in other words up to a homotopy of complexes. 
Thus, for any choice of $F'$ which is a complex of bundles representing $F$ up to quasiisomorphism, we get
an action of $\cL$. This will induce the same map as the original $\cL \otimes F\rightarrow F$ on cohomology,
and the same is true after applying a derived functor such as taking the
fiber over a point $x\in J$.
\end{proof}

Replace the notation $F'$ above by the shorter notation $F$ in what follows. See \cite{Cailotto} for a
similar statement with application to a Steenbrink-type theorem. 

\begin{proposition}
\label{criterion}
Suppose $F$ is a complex of vector bundles over an affine curve or a complex disk $S$. Let $x\in S$ 
and let $J$ be the divisor $[x]$. Define $\cL$ as before with respect to this divisor. 
Suppose that we are given an action 
$$
M:\cL \otimes _{\cO _S}F\rightarrow F
$$
where the tensor product uses the right module structure on $\cL$ and leaves the left module structure in the answer.
Suppose that, when restricted to $\cO _S\subset \cL$ it gives a map
$$
M_0 : F = \cO _S \otimes _{\cO _S} F \rightarrow F
$$
which induces the identity on cohomology sheaves and on the cohomology spaces of $F_x$. 
Suppose furthermore that for our point $x$ the morphism
$$
M_x : \cL _x\otimes _{\cO _S}F_x\rightarrow F_x
$$
induces a nilpotent action of $( t\frac{\partial}{\partial t} ) _x \subset \cL _x$ on 
$F_x$. Then, locally in a neighborhood of $x$ the complex $F$ splits as a direct sum of its cohomology sheaves,
and these cohomology sheaves are locally free. 
\end{proposition}
\begin{proof}
The map $M$ gives a map
$$
H^i(M) : \cL \otimes _{\cO _S}\cH ^i(F)\rightarrow \cH ^i(F)
$$
on cohomology sheaves. Away from $D$, this map is a connection,
and as is well-known it implies that the cohomology sheaves are locally free. Thus our only problem is to prove that
they are locally free near a point $x\in J$.

Over the disc or affine curve $S$  the complex $F$ is quasiisomorphic to the direct sum of its cohomology sheaves.
This is because there are no $Ext ^i$ terms for $i\geq 2$ over $\cO _S$ since $S$ is one-dimensional.

We can then choose a minimal resolution of each cohomology sheaf. This can be done explicitly, because by the
Chinese remainder theorem we have that $H^i(F)$ is a direct sum of locally free modules and
modules of the form $\cO _S / (z^n)$. For the latter pieces if they exist (our goal is to show that
they don't occur) choose as resolution
$$
0 \rightarrow \cO _S \stackrel{z^n}{\rightarrow} \cO _S\rightarrow \cO _S / (z^n)\rightarrow 0.
$$
In this way we obtain a quasiisomorphism $R \rightarrow F \rightarrow R$ between $F$ and a complex $R$ of locally
free sheaves with the ``minimality'' property that the differentials in $R$ vanish at our singular point $x\in D$.

Now replace $F$ by the minimal resolution $R$. We still get a map
$$
M^R:\cL \otimes _{\cO _S} R \rightarrow R.
$$
Its restriction $M^R_0$ to $\cO _S\subset \cL$ is a map $M_0:R\rightarrow R$ which induces the identity
on cohomology sheaves. The fiber $R_x$ is a complex whose differentials are equal to zero, in particular
$$
R_x \cong \bigoplus _i H^i(F_x).
$$
Therefore, $M^R_{0,x}$ is equal to the identity by our hypothesis that $M_0$ induces the identity
on $H^i(F_x)$. In particular, $M^R_{0,x}$  is invertible, so
by Nakayama's lemma, after possibly restricting the size of our disc or affine neighborhood of $x$, 
we may assume that $M^R_0 : R\rightarrow R$ is
invertible. Define
$$
M':= (M^R_0)^{-1}\circ M^R : \cL \otimes _{\cO _S} R \rightarrow R.
$$
We have now succeeded in normalizing so that the restrition of $M'$ to $\cO _S\subset \cL$ is the identity.
One can see, using \eqref{eq.-c}, that this exactly defines a logarithmic connection
denoted $\nabla$ on
the whole complex $R$. In particular each $R^i$ now has a logarithmic connection and the differentials of $R$ are
compatible with these connections.

Note that the fiber $\cL _x$ is a vector space which is well-defined with respect to either the left or right
$\cO _S$-module structures of $\cL$ because the symbols $\sigma (v)$ of elements of $\cL$ vanish at $x\in D$.
The exact sequence for $\cL$ splits over $x$ to give just
$$
\cL _x \cong \comx \oplus T_S(-x)_x \cong \comx \oplus \comx.
$$
Thus $M$ gives a well-defined map
$$
M_x (t\frac{\partial}{\partial t})_x:  F_x \rightarrow F_x.
$$
The hypothesis of the present proposition is that this endomorphism of $F_x$ induces nilpotent endomorphisms
on cohomology.  
Composing with $R_x\rightarrow F_x \rightarrow R_x$ which are again quasiisomorphisms,
gives the map
$$
M^R_x (t\frac{\partial}{\partial t})_x:  R_x \rightarrow R_x.
$$
However, $R_x \cong H^{\cdot} (F_x)$. Thus the action of $M^R_x (t\frac{\partial}{\partial t})_x$ on
$R_x$ is nilpotent. 

We claim that
the connection $\nabla $ on the complex $R$ has nilpotent residues at $x$.
Indeed, the residue of $\nabla$ at $x$ is equal to $(M^R_{0,x})^{-1}\circ M^R_x (t\frac{\partial}{\partial t})_x$,and 
we know from the hypothesis of our proposition that 
$M^R_{0,x}$ is the identity, so the residue of $\nabla$ is nilpotent. 

The claim implies, in particular, that each $R^i$ is a Deligne canonical extension.

\begin{lemma}
The differentials of the complex $R$ are zero.
\end{lemma}
\begin{proof}
Look at any one of the differentials
$$
d_i : R^i\rightarrow R^{i+1}.
$$
This map is a map of vector bundles with logarithmic connection. Thus $d_i$ may be considered as a flat section
of the vector bundle with logarithmic connection $Hom (R^i, R^{i+1})$. This latter is also a Deligne canonical extension.
On the other hand, $d_i (x)=0$ by the minimality assumption
for $R$. We are now in the following situation: we have a vector bundle with logarithmic connection which is a
Deligne canonical extension (i.e. the residue is nilpotent), and we have a flat section which vanishes at $x$.
We claim that this implies that the section vanishes everywhere. The claim is clearly true when we are
dealing with sections of the canonical extension of a trivial connection. However, as we are working over a disc,
our Deligne extension is a successive extension of trivial pieces, so by an induction argument, we get that $d_i=0$.
\end{proof}

\begin{corollary}
The complex $F$ is locally a direct sum of its cohomology sheaves, and these cohomology sheaves are locally free.
\end{corollary}
\begin{proof}
This statement is invariant under quasiisomorphism, and it is true for $R$ because of the previous corollary
so it is true for $F$.
\end{proof}

This corollary completes the proof of the proposition. 
\end{proof}

{\em Proof of Theorem \ref{th.-steen}:}
Lemma \ref{actionL} constructs the Gauss-Manin connection on the cohomology complex
and shows that it satisfies the condition that the residual actions over points
$x\in J$ are nilpotent, and in Proposition \ref{criterion} we have seen that this nilpotence property
implies that the cohomology complex $F$ is locally a direct sum of locally free cohomology sheaves.
Along the way, the nilpotence property provides the statement of the second part of Theorem \ref{th.-steen}. 
This completes the proof when $S$ is a curve, and as stated at the beginning of the section, it is 
enough to treat the case when the base is a curve.  The last statement
about Chern characters follows because a complex and the direct sum of its cohomology are equivalent in $K$-theory.
\eop


\begin{thebibliography}{AAAA}

\bibitem[Ab-dJ]{Ab-dJ}Abramovich, D., de Jong, A.J. {\em Smoothness, semistability, and toroidal geometry}, 
J.Algebraic Geom. \textbf{6} (1997), 789-801.

\bibitem[AGV]{AbramovichGraberVistoli} Abramovich, D., Graber, T., Vistoli, A.  
{\em Algebraic orbifold quantum products}, 
Orbifolds in mathematics and physics (Madison, WI, 2001), 1--24, Contemp. Math., 
\textbf{310}, Amer. Math. Soc., Providence, RI, 2002.

\bibitem[AK]{AbramovichKaru}
Abramovich, D., Karu, K.
{\em Weak semistable reduction in characteristic $0$}.
Invent. Math. \textbf{139} (2000), 241--273.

\bibitem[Be]{Berthelot} Berthelot, P.  
{\em Cohomologie cristalline des sch\'emas de caract\'eristique $p>0$.}  
{\sc Lect. Notes in Maths.} {\bf 407},
Springer-Verlag, Berlin-New York (1974). 

\bibitem[Bi]{Biswas} Biswas, I. {\em Chern classes for parabolic bundles}, 
J. Math. Kyoto Univ.  37  (1997),  no. \textbf{4}, 597--613.

\bibitem[Bi2]{Biswas2} Biswas, I. 
{\em Parabolic bundles as orbifold bundles}, Duke Math. J., {\bf 88} (1997), 305-325.


\bibitem[Bi-Iy]{Bi-Iy} Biswas, I., Iyer, J. {\em Vanishing of Chern classes of the de Rham bundles for
some families of moduli spaces}, preprint 2004.

\bibitem[BE]{BE} Bloch, S., Esnault, H.  {\em Algebraic Chern-Simons
theory}, Amer. J. Math. 119 (1997),  no. \textbf{4}, 903--952.

\bibitem[BE2]{BE2}  Bloch, S., Esnault, H. 
{\em A Riemann-Roch theorem for flat bundles, with values in the algebraic Chern-Simons theory}, 
Ann. of Math. (2) 151 (2000), no. \textbf{3}, 1025--1070. 

\bibitem[Bo]{Borne} Borne, N. {\em Fibr\'es paraboliques et champ des racines}, preprint 2005, more recently 
posted as {\tt math.AG/0604458}.

\bibitem[Cad]{Cadman} Cadman, C. {\em Using stacks to impose tangency conditions on curves}, Preprint {\tt math.AG/0312349}.

\bibitem[Cai]{Cailotto}
Cailotto, M.
{\em Algebraic connections on logarithmic schemes.}
C. R. Acad. Sci. Paris S\'er. I Math. \textbf{333} (2001),  1089--1094.

\bibitem[Ch-Ru]{ChenRuan}  Chen, W., Ruan, Y.  {\em Orbifold Gromov-Witten theory}, Orbifolds in mathematics and physics 
(Madison, WI, 2001), 25--85, Contemp. Math., \textbf{310}, Amer. Math. Soc., Providence, RI, 2002.

\bibitem[CK]{ConsaniKim}
Consani, C., Kim, M.
{\em The Picard-Lefschetz formula and a conjecture of Kato: the case of Lefschetz fibrations.}
Math. Res. Lett. \textbf{9} (2002), 621--631.

\bibitem[dJ]{deJong} de Jong, A.J.
{\em Smoothness, semi-stability and alterations}, 
{Publ. Math. I. H. E. S.} {\bf 83} (1996), 51-93.

\bibitem[De1]{De} Deligne, P. {\em Equations diff\'erentielles a points singuliers
reguliers}. {\sc Lect. Notes in Math.} $\bf{163}$, 1970.

\bibitem[De2]{De2} Deligne, P. {\em Th\'eor\`eme de Lefschetz et crit\`eres de d\'eg\'en\'erescence de suites
 spectrales}, Publ.Math.IHES \textbf{35} (1968), 107-126.
 
 
\bibitem[De3]{De3} Deligne, P. {\em Th\'eorie de Hodge. III}, Inst. Hautes \'Etudes Sci. 
Publ. Math. No. \textbf{44} (1974), 5--77. 

\bibitem[El]{Elzein} Elzein, F. {\em Hodge-de Rham theory with degenerating coefficients}, 
Preprint \verb}math.AG/0311083}.

\bibitem[Es1]{Es1} Esnault, H.  {\em Recent developments on characteristic
classes of flat bundles on complex algebraic manifolds}, Jahresber. Deutsch. Math.-Verein.  
98  (1996),  no. $\bf{4}$, 182--191.

\bibitem[Es2]{Es2} Esnault, H.  {\em Private communication.}

\bibitem[Es-Vi1]{Es-Vi1} Esnault, H., Viehweg, E.  {\em Logarithmic De Rham complexes
and vanishing theorems}, Invent.Math., $\bf{86}$, 161-194, 1986.

\bibitem[Es-Vi2]{Es-Vi2} Esnault, H., Viehweg, E.  {\em Chern classes of Gauss-Manin bundles of weight 1 vanish}, 
$K$-Theory  26  (2002),  no. \textbf{3}, 287--305.

\bibitem[Fa]{Faltings} Faltings, G. {\em $F$-isocrystals on open varieties: results and conjectures},
{\em The Grothendieck Festschrift, Vol. II} 219--248,
{\sc Progr. Math.} {\bf 87} Birkh\"auser Boston, Boston, MA, 1990. 

\bibitem[Fa2]{Faltings2} Faltings, G. 
{\em Crystalline cohomology and $p$-adic Galois-representations},
Algebraic analysis, geometry, and number theory (Baltimore, MD, 1988), 25--80,
Johns Hopkins Univ. Press, Baltimore, MD, 1989. 

\bibitem [Fu]{Fu} Fulton, W. {\em Intersection theory}, Second edition. Ergebnisse der Mathematik und ihrer Grenzgebiete. 3. Folge, \textbf{2}. Springer-Verlag, Berlin, 1998. xiv+470 pp.

\bibitem[Gi]{Gi} Gillet, H. {\em  Intersection theory on algebraic stacks and $Q$-varieties}, Proceedings of the 
Luminy conference on algebraic $K$-theory (Luminy, 1983).
J. Pure Appl. Algebra \textbf{34} (1984), 193--240.

\bibitem[Gr]{GrothDR} Grothendieck, A. {\em On the de Rham cohomology of algebraic varieties},
{\em Publ. Math. I. H.E.S.} {\bf 29} (1966), 95-103.

\bibitem[Ha]{HartshorneManuscripta} Hartshorne, R.
{\em Algebraic de Rham cohomology}, Manuscripta Math., {\bf 7} (1972), 125-140.

\bibitem[Ha2]{HartshorneIHES} Hartshorne, R.
{\em On the De Rham cohomology of algebraic varieties},
{\em Publ. Math. I.H.E.S.} {\bf 45} (1975), 5-99.

\bibitem[HL]{HerreraLieberman}
Herrera, M.; Lieberman, D.
{\em Duality and the de Rham cohomology of infinitesimal neighborhoods},
Invent. Math. \textbf{13} (1971), 97--124.

\bibitem[Hi]{Hi}  Hironaka, H. {\em Resolution of singularities of an algebraic variety over a field 
of characteristic zero. I, II.}, Ann. of Math. (2) \textbf{79} (1964), 109--203.

\bibitem[Il]{Illusie} Illusie, L. 
{\em R\'eduction semi-stable et d\'ecomposition de complexes de de Rham \`a coefficients},  Duke Math. J., 
{\bf 60} (1990), 139-185.


\bibitem[IKN]{IllusieKatoNakayama}
Illusie, L., Kato, K., Nakayama, C.
{\em Quasi-unipotent logarithmic Riemann-Hilbert correspondences}, J. Math. Sci. Univ. Tokyo, {\bf 12} (2005), 1-66.

\bibitem[In]{Inaba} Inaba, M. {\em Moduli of parabolic connections on a curve and Riemann-Hilbert correspondence},
Preprint \verb}math.AG/0602004}.

\bibitem[IIS]{InabaEtAl} Inaba , M.(Kyushu), Iwasaki, K., Saito, Masa-Hiko.
{\em Moduli of Stable Parabolic Connections, Riemann-Hilbert correspondence and Geometry of 
Painlev\'{e} equation of type VI, Part I},
Preprint \verb}math.AG/0309342}.

\bibitem[Iy]{Iy} Iyer, J.  {\em The de Rham bundle on a compactification
of moduli space of abelian varieties.}, Compositio Math.  136  (2003),  no. $\bf{3}$, 317--321.

\bibitem[Ko]{FKato} Kato, F.
{\em The relative log Poincar\'e lemma and relative log de Rham theory}, Duke Math. J., {\bf 93} (1998), 179-206.

\bibitem[KMN]{KatoMatsubaraNakayama}
Kato, K., Matsubara, T., Nakayama, C.
{\em Log $C^{ \infty}$-functions and degenerations of Hodge structures.}
Algebraic geometry 2000, Azumino (Hotaka), 
{\sc Adv. Stud. Pure Math.}, \textbf{36},
Math. Soc. Japan, Tokyo (2002), 269--320

\bibitem[KN]{KatoNakayama}
Kato, K., Nakayama, C.
{\em Log Betti cohomology, log \'etale cohomology, and log de Rham cohomology of log schemes over $\comx$}.
Kodai Math. J. \textbf{22} (1999), 161--186.

\bibitem[Kz]{Ka} Katz, N. M. {\em \'Etude cohomologique des pinceaux de Lefschetz}, SGA 7, expos\'e 18, LNM
\textbf{340}, Springer.

\bibitem [Kz2]{KatzICM} Katz, N. M. {\em The regularity theorem in algebraic geometry}, 
Actes du Congr\'es International des Math\'ematiciens (Nice, 1970), Tome \textbf{1}, pp. 437--443. 
Gauthier-Villars, Paris, 1971. 

\bibitem[KO]{KatzOda}
Katz, N. M., Oda, T.
{\em On the differentiation of de Rham cohomology classes with respect to parameters}.
J. Math. Kyoto Univ. {\bf 8}, 2 (1968), 199-213.
 
\bibitem[Li]{Li} Li, Jiayu.
{\em Hermitian-Einstein metrics and Chern number inequalities on parabolic stable bundles over K\"ahler manifolds}
Comm. Anal. Geom. 8 (2000), no. \textbf{3}, 445--475.

\bibitem[Li-Na]{LiNarasimhan} Li, Jiayu., Narasimhan, M. S.
{\em Hermitian-Einstein metrics on parabolic stable bundles}
Acta Math. Sin. (Engl. Ser.) {\bf 15} (1999),  93-114.

\bibitem[Man]{Manin} Manin, Y. {\em  Moduli Fuchsiani}, Annali Scuola Normale Sup. di Pisa Ser. III
{\bf 19} (1965), 113-126.

\bibitem[Ma-Ol]{MatsukiOlsson} Matsuki, K., Olsson, M. {\em  Kawamata-Viehweg 
vanishing and Kodaira vanishing for stacks}, Math. Res. Lett.
{\bf 12} (2005), 207-217.

\bibitem[Ma-Yo]{MaruyamaYokogawa} Maruyama, M., Yokogawa, K. 
{\em Moduli of parabolic stable sheaves},
Math. Ann. 293 (1992), no. \textbf{1}, 77--99.

\bibitem[Mo]{Mochizuki} Mochizuki, T. {\em Kobayashi-Hitchin correspondence for 
tame harmonic bundles and an application}, 
Preprint {\tt math.DG/0411300}.

\bibitem[Mu]{Mu} 
Mumford, D. {\em Towards an Enumerative Geometry of the
Moduli Space of Curves}, Arithmetic and geometry, Vol. II, 271--328, Progr.
Math., $\bf{36}$, Birkh$\ddot{a}$user Boston, Boston, MA, 1983.

\bibitem[Og]{Ogus}
A. Ogus. 
{\em $F$-crystals, Griffiths transversality, and the Hodge decomposition.}
Ast\'erisque \textbf{221} (1994).

\bibitem[Oh]{Oh} Ohtsuki, M. {\em A residue formula for Chern classes associated with logarithmic connections}, Tokyo J. Math. 5 (1982), no. \textbf{1}, 13--21.  

\bibitem[Pa]{Panov} Panov, D.  {\em Doctoral thesis}, 2005.

\bibitem[Sa]{Saito} Saito, M.
{\em Mixed Hodge modules}, Proc. Japan Acad. Ser. A Math. Sci., {\bf 62} (1986), 360-363.

\bibitem[Se]{Seshadri} 
Seshadri, C. S.
{\em Moduli of vector bundles on curves with parabolic structures.}
Bull. Amer. Math. Soc. \textbf{83} (1977), 124--126.

\bibitem[Sr]{Sr} Srinivas, V. {\em Algebraic $K$-theory}, Progress in Mathematics, \textbf{90}. 
Birkhauser Boston, Inc., Boston, MA, 1996. xviii+341 pp.

\bibitem[St]{Steen} Steenbrink, J.  {\em Limits of Hodge structures}, 
Invent. Math.  31  (1975/76), no. \textbf{3}, 229--257.


\bibitem[St-Wr]{SteerWren} Steer, B., Wren, A.
{\em The Donaldson-Hitchin-Kobayashi correspondence for parabolic bundles over orbifold surfaces}, 
Canad. J. Math. 53 (2001), no. \textbf{6}, 1309--1339.

\bibitem[To]{Toen} Toen, B.  {\em Th\'eor\`emes de Riemann-Roch pour les champs de Deligne-Mumford}  
$K$-Theory 18 (1999), no. \textbf{1}, 33--76.


\bibitem[Ts]{Tsuji} Tsuji, T.
{\em On $p$-adic nearby cycles of log smooth families}, Bull. Soc. Math. France, {\bf 128} (2000), 529-575.

\bibitem[Ts2]{Tsuji2} Tsuji, T.
{\em $p$-adic \'etale cohomology and crystalline cohomology in the semi-stable reduction case}, 
Invent. Math., {\bf 137} (1999), 233-411.


\bibitem [vdG]{vdG} van der Geer, G.  {\em Cycles on the moduli space of abelian
varieties}, Moduli of curves and abelian varieties, 65-89, Aspects Math. E33,
 Vieweg, Braunschweig 1999.


\bibitem[Vi]{Vi} Vistoli, A. {\em Intersection theory on algebraic stacks and on their moduli spaces}, 
Invent. Math. 97 (1989), no. \textbf{3}, 613--670.

\bibitem[Yo]{Yokogawa} Yokogawa, K.
{\em Compactification of moduli of parabolic sheaves and moduli of parabolic Higgs sheaves},
J. Math. Kyoto Univ. 33 (1993), no. \textbf{2}, 451--504.

\end {thebibliography}
\end{document}